\documentclass[letter,11pt]{article}

\usepackage[english]{babel}
\usepackage[utf8]{inputenc}

\usepackage[margin=1.2 in, top=1.1 in, bottom= 1.2 in]{geometry}
\makeatletter
\g@addto@macro\normalsize{%
  \setlength\abovedisplayskip{7pt}
  \setlength\belowdisplayskip{7pt}
  \setlength\abovedisplayshortskip{7pt}
  \setlength\belowdisplayshortskip{7pt}
}
\makeatother

\usepackage{tocloft}

\usepackage{amsfonts}
\usepackage{mathrsfs}
\usepackage{bbm}
\usepackage{latexsym}
\usepackage{math dots}
\usepackage{amssymb}
\usepackage{mathtools}
\usepackage{relsize}
\usepackage{lipsum}

\usepackage{todonotes}
\usepackage{enumitem}
\setlist{nolistsep} 	%[label=(\roman{enumi}),ref=(\roman{enumi}),leftmargin=*]
\usepackage{amsthm}

\usepackage{xcolor}
\definecolor{Color1}{rgb}{0.0, 0.42, 0.47}%Blue
\definecolor{Color2}{rgb}{0.78, 0.11, 0.0}%Scarlet

\usepackage{titlesec}
\titleformat{\section}
  {\Large\center\bfseries}
  {\thesection.}{.7em}{}
\titlespacing*{\section}{0pt}{3.5ex plus 0ex minus 0ex}{1.5ex plus 0ex}
\titleformat{\subsection}
  {\center\bfseries}
  {\thesubsection.}{.7em}{}
\titlespacing*{\subsection}{0pt}{3.5ex plus 0ex minus 0ex}{1.5ex plus 0ex}
\titleformat{\subsubsection}
  {\center\bfseries}
  {\thesubsubsection.}{.7em}{}
\titlespacing*{\subsubsection}{0pt}{3.5ex plus 0ex minus 0ex}{1.5ex plus 0ex}

\addto\captionsenglish{}

\usepackage{titling}
\usepackage[linktocpage=true]{hyperref}
\usepackage[capitalize]{cleveref}
\hypersetup{citecolor = Color1,colorlinks,
			linkcolor = black,
			urlcolor = Color2}

\crefalias{maintheorem}{theorem}

\newtheoremstyle{plain}{3mm}{3mm}{\slshape}{}{\bfseries}{.}{.5em}{}
\newtheoremstyle{definition}{2mm}{2mm}{}{}{\bfseries}{.}{.5em}{}
\theoremstyle{plain} %PLAIN
	
\newtheorem{theorem}{Theorem}[section]

\newtheorem{proposition}[theorem]{Proposition}

\newtheorem{conjecture}[theorem]{Conjecture}
\newtheorem{lemma}[theorem]{Lemma}
\newtheorem{corollary}[theorem]{Corollary}
\theoremstyle{definition} %DEFINITION
\newtheorem{defn}[theorem]{Definition}
\newtheorem{remark}[theorem]{Remark}

\theoremstyle{plain} 
\newcounter{MainTheoremCounter}

\newtheorem{maintheorem}[MainTheoremCounter]{Theorem}

\theoremstyle{plain} 
\newcounter{MainConjectureCounter}

\theoremstyle{plain}
\newtheorem*{namedthm}{\namedthmname}
\newcounter{namedthm}
	
\usepackage{chngcntr}
\numberwithin{equation}{section}

\allowdisplaybreaks

\newcommand{\Cesaro}{Ces\`{a}ro}
\newcommand{\Erdos}{Erd\H{o}s}

\newcommand{\Turan}{Tur{\'a}n}
\newcommand{\Matomaki}{Matom\"{a}ki}
\newcommand{\Radziwill}{Radziwi\l\l{}}
\newcommand{\Halasz}{Hal\'{a}sz}

\newcommand{\Oh}{{\mathrm O}}
\newcommand{\oh}{{\mathrm o}}
\newcommand{\N}{\mathbb{N}}
\newcommand{\Z}{\mathbb{Z}}
\newcommand{\R}{\mathbb{R}}
\newcommand{\C}{\mathbb{C}}
\newcommand{\Q}{\mathbb{Q}}
\newcommand{\T}{\mathbb{T}}

\renewcommand{\epsilon}{\varepsilon}
\renewcommand{\leq}{\leqslant}
\renewcommand{\geq}{\geqslant}
\renewcommand{\setminus}{\backslash}
\renewcommand{\Re}{\mathrm{Re}}
\renewcommand{\Im}{\mathrm{Im}}
\renewcommand{\P}{\mathbb{P}}
\renewcommand{\subset}{\subseteq}

\renewcommand{\d}{~\mathrm{d}}

\newcommand{\lio}{\boldsymbol{\lambda}}

\usepackage[normalem]{ulem}

\usepackage[makeroom]{cancel}

\newcommand{\E}{\raisebox{-0.4ex}{\scalebox{1.4}[1.3]{$\mathbb{E}$}}}
\newcommand{\Esub}[1]{\raisebox{-0.4ex}{\scalebox{1.4}[1.3]{$\mathbb{E}$}}\raisebox{-0.35ex}{\ensuremath{_{\scriptstyle #1}}}}
\newcommand{\logE}{{\E}^{\log}}

\newcommand{\Dist}{\mathbb{D}}
\newcommand{\1}{\mathbbm{1}}

\newcommand{\f}{f_{\xi,N}}
\newcommand{\g}{g_{\xi,N}}
\newcommand{\Pc}{\mathcal{P}}

\usepackage[nottoc]{tocbibind}

\author{By~~{\scshape Dimitrios~Charamaras}~~and~~{\scshape Florian~K.~Richter}}
\date{\small \today}
\title{\bfseries 
Asymptotic independence of $\Omega(n)$ and $\Omega(n+1)$ along logarithmic averages}

\begin{document}

\maketitle

\begin{abstract}
Let $\Omega(n)$ denote the number of prime factors of a positive integer $n$ counted with multiplicities. We show that for any bounded functions $a,b\colon\N\to\C$,
$$\frac{1}{\log{N}}\sum_{n=1}^N \frac{a(\Omega(n))b(\Omega(n+1))}{n}
= \Bigg(\frac{1}{N}\sum_{n=1}^N a(\Omega(n))\Bigg)\Bigg(\frac{1}{N}\sum_{n=1}^N b(\Omega(n))\Bigg) + \oh_{N\to\infty}(1).$$
This generalizes a theorem of Tao on the logarithmically averaged two-point correlation Chowla conjecture.
Our result is quantitative and the explicit error term that we obtain establishes double-logarithmic savings.
As an application, we obtain new results about the distribution of $\Omega(p+1)$ as $p$ ranges over $\ell$-almost primes for a ``typical'' value of $\ell$.
\end{abstract}

\tableofcontents
\thispagestyle{empty}

%==========================================================
%==========================================================

%SECTION
\section{Introduction}\label{Intro}

For a positive integer $n\in\N=\{1,2,3,\ldots\}$, let $\Omega(n)$ denote the number of prime factors of~$n$ counted with multiplicities.
Given $N\in\N$ we use the notation $[N]$ for the set of integers $\{1,\ldots,N\}$, we denote the \Cesaro{} average of a function $f\colon\N\to\C$ over a finite non-empty set $A\subset\N$ as
$$\Esub{n\in A}f(n) = \frac{1}{|A|}\sum_{n\in A}f(n),$$
and the logarithmic average of $f$ over $A$ as 
$$\logE_{n\in A}f(n)=\frac{1}{\sum_{n\in A}\frac{1}{n}}\sum_{n\in A} \frac{f(n)}{n}.$$

The starting point of our discussion is an old conjecture of Chowla concerning the Liouville function $\lio(n)\coloneqq (-1)^{\Omega(n)}$.
It asserts that asymptotically there is no correlation among the sign patterns of $\lio(n)$ along consecutive integers.

\begin{conjecture}[Chowla {\cite{Chowla87}; see also \cite{Ramare18}}]
For any $k\in\N$,
\begin{equation}
\label{eqn_chowla}
    \Esub{n\in[N]} \lio(n) \lio(n+1)\cdots \lio(n+k-1) = \oh_{N\to\infty}(1).
\end{equation}
\end{conjecture}

If the \Cesaro{} average in \eqref{eqn_chowla} is replaced by a logarithmic average, then the case $k=2$ of \eqref{eqn_chowla} was proved in \cite{Tao16}.

\begin{theorem}[{cf.~\cite[Theorem 2]{Tao16}}]\label{thm_log_2_correlation_chowla}
We have
\[
\logE_{n\in[N]} \lio(n)\lio(n+1)= \oh_{N\to\infty}(1).
\]
\end{theorem}

In \cite{helf_radz}, Helfgott and \Radziwill{} established double-logarithmic savings in the above, while the best known result to date is due to Pilatte \cite{pilatte2023improved} who proved \cref{thm_log_2_correlation_chowla} with logarithmic savings.

The method used in \cite{Tao16} was later extended in \cite{TT18} to establish the logarithmically averaged Chowla's conjecture for any odd $k$. 
Albeit significant recent progress combining ideas from number theory, higher order Fourier analysis, and ergodic theory (cf.\ \cite{Sarnak11,EAKLdlR17,Frantzikinakis17,tao2017equivalence,FH18a,FKL18,GLdlR21,KMT23_almost_all_scales,MRTT2023_higher_uniformity,MSTT2023,KLRT2024,Matomaki-Radziwill-Shao-Tao-Teravainen2026}),
Chowla's conjecture, and the more general Elliott's conjecture\footnote{The original Elliott's conjecture stated by Elliott in \cite{Elliott92} is false (see \cite[Theorem B.1]{MRT15} for the counterexample). The corrected version of the conjecture was given in \cite{MRT15} and it is still open.} (cf.~\cite{Elliott92,MRT15}), remain open for other $k$, even in the case of logarithmic averages.

A different relation for the sequence $\Omega(n)$, introduced in \cite{Richter_PNT_21}, draws inspiration from an ergodic-theoretic perspective.

\begin{theorem}[cf.~{\cite{Richter_PNT_21}}]
\label{thm_Omega_is_good_for_UE}
For any bounded arithmetic function $a\colon\N\to\C$, we have
\begin{equation}
\label{eqn_pnt_1}
\Esub{n\in[N]} a(\Omega(n)) = \Esub{n\in[N]} a(\Omega(n)+1) +\oh_{N\to\infty}(1).
\end{equation}
\end{theorem}

We remark that if $a(n)=(-1)^n$ then \eqref{eqn_pnt_1} reduces to the $k=1$ case of Chowla's conjecture, %coincides with the classical Liouville function $\lio(n)=(-1)^{\Omega(n)}$ and \eqref{eqn_pnt_1} is equivalent to $\frac{1}{N}\sum_{n=1}^N\lio(n) =\oh_{N\to\infty}(1)$,
which is a well-known equivalent form of the Prime Number Theorem.
For generalizations of \eqref{eqn_pnt_1} and its connections to ergodic theory see~\cite{BR22,Loyd23, Charamaras_MVT_24}. 
Theorems~\ref{thm_log_2_correlation_chowla}~and~\ref{thm_Omega_is_good_for_UE} highlight different number-theoretic aspects of the distribution of the values of $\Omega(n)$ and 
our goal is to unify these approaches.
% and opening up unexplored directions of research.
Instead of restricting to the joint distribution of multiplicative parity values, as in Chowla's conjecture, it is natural to inquire about more general forms of independence regarding the joint distribution of $\Omega(n),\Omega(n+1),\ldots,\Omega(n+k-1)$. 
Based on \cref{thm_Omega_is_good_for_UE}, we propose the following conjecture.

\begin{conjecture}
\label{conj_functional_chowla}
For any $k\in\N$ and any bounded functions $a_1,\ldots,a_k\colon\N\to\C$,
$$
\Esub{n\in[N]} a_1(\Omega(n))a_2(\Omega(n+1))\cdots a_k(\Omega(n+k-1))
= \prod_{i=1}^k \Esub{n_i\in[N]} a_i(\Omega(n_i)) + \oh_{N\to\infty}(1).
$$
\end{conjecture}

Notice that if $a_1(n)=\ldots=a_k(n)=(-1)^{n}$, then \cref{conj_functional_chowla} reduces to Chowla's conjecture. More generally, if $a_1,\ldots,a_k\colon\N\to\mathbb{S}^1=\{z\in\C: |z|=1\}$ are group characters of $(\N,+)$, then \cref{conj_functional_chowla} becomes a special case of Elliott's conjecture. The intended upgrade of \cref{conj_functional_chowla} is that no structural restrictions, particularly no algebraic ones, are imposed on the functions $a_1,\ldots,a_k$.

\subsection{Main results}

The aim of this paper is to validate the logarithmically averaged two-point correlation case of \cref{conj_functional_chowla}. %\new{For simplicity, we focus on the case $(h_1,h_2)=(0,1)$.} 
Our main result is as follows:

\begin{maintheorem}\label{thm_A}
For any bounded functions $a,b\colon\N\to\C$, we have
    \begin{equation}\label{thm_A_eq}
        \logE_{n\in[N]} a(\Omega(n))b(\Omega(n+1)) 
        = \Big(\Esub{n\in[N]}a(\Omega(n))\Big)\Big(\Esub{n\in[N]}b(\Omega(n))\Big)
        + \Oh\bigg(\frac{(\log\log\log N)^2}{\sqrt{\log\log N}}\bigg).
    \end{equation}
    Moreover, the error term only depends on $\max\{\|a\|_\infty,\|b\|_\infty\}$.
\end{maintheorem}

Theorem~\ref{thm_A} is a generalization of \cref{thm_log_2_correlation_chowla}, since choosing $a(n)=b(n)=(-1)^n$ in \eqref{thm_A_eq} implies \eqref{eqn_chowla}.

\begin{remark}
We remark that the optimal error term in \eqref{thm_A_eq} is 
\[
\Oh\bigg(\frac{1}{\sqrt{\log\log{N}}}\bigg),
\]
since this bound is known to be best possible for the \Erdos{}-Kac theorem due to the work of \cite{renyi-turan}, and \eqref{thm_A_eq} contains the (logarithmically averaged) \Erdos{}-Kac theorem as a special case. Consequently, the error term in Theorem~\ref{thm_A} approaches the optimal bound, falling short by a factor of only $(\log\log\log N)^2$.
\end{remark}

\begin{remark}
After the first version of this article appeared on arXiv, Harald A.~Helfgott pointed out to us that one can obtain a different proof of Theorem~\ref{thm_A} by replacing our analysis of the small and large frequencies (see \cref{proof strategy} for details) with \cite[Proposition 8.12]{helf_radz}. This alternative approach has the advantage that it yields the optimal error term $\Oh\Big(\frac{1}{\sqrt{\log\log{N}}}\Big)$ instead of $\Oh\Big(\frac{(\log\log\log N)^2}{\sqrt{\log\log{N}}}\Big)$ obtained via our method. Since our original proof was developed independently, we have retained it in the main text; for completeness, we present the alternative argument in \cref{appendix}.
\end{remark}

Theorem~\ref{thm_A} demonstrates a natural form of asymptotic independence between $\Omega(n)$ and $\Omega(n+1)$.
We can provide a clearer, more intuitive, understanding of this behavior by recasting our result in the language of probability theory, in the spirit of the \Erdos{}-Kac theorem.
Recall that if $X$ is a Gaussian random variable with mean $\mu=\mathbb{E}[X]$ and standard deviation $\sigma=\sqrt{\mathbb{E}[(X-\mu)^2]}$ then the probability density function of $X$ is given by
\begin{equation}
\label{eqn_Gaussian_prob_dens_fctn}
f(x;\mu,\sigma)= \frac{1}{\sigma\sqrt{2\pi}}\,e^{-\frac{1}{2}\left(\frac{x-\mu}{\sigma}\right)^2}.
\end{equation}
This means that for any bounded and continuous $g\colon \R\to \C$ we have
\begin{equation}
\label{eqn_Gaussian_density_function_characterization}
\mathbb{E}[g(X)] = \int_{-\infty}^\infty g(x)f(x;\mu,\sigma) \d x. 
\end{equation}
In analogy, the ``mean'' and ``standard deviation'' of $\Omega(n)$ over $[N]$ are defined as
\begin{equation}\label{mu_N and sigma_N}
    \mu_N = \log\log{N} \sim \Esub{n\in[N]}\Omega(n)
    \quad\text{and}\quad
     \sigma_N = \sqrt{\log\log{N}} \sim \Big(\Esub{n\in[N]}\big|\Omega(n)-\mu_N\big|^2\Big)^\frac{1}{2}.
\end{equation}

Given $A>1$, we define the {\em typical range} $T_{A,N}$ of $\Omega(n)$ as
\begin{equation}\label{typical_range}
    T_{A,N} = [\mu_N-A\sigma_N,\mu_N+A\sigma_N] = [\log\log{N}-A\sqrt{\log\log{N}}, \log\log{N}+A\sqrt{\log\log{N}}].
\end{equation}
Let $\P_\ell=\{n\in\N: \Omega(n)=\ell\}$ denote the set of all $\ell$-almost primes.
Based on the uniform asymptotic estimates of $\P_\ell\cap[N]$ obtained by \Erdos{} in \cite{Erdos48a} (see also \cite{Selberge_Sathe}), one can show (see \cite[Lemma 3.4]{Loyd23} and \cref{pi_k_estimate} below) that for all $A>1$ and uniformly over all $\ell$ in the typical range $T_{A,N}$, one has
\begin{equation*}
\label{eqn_loy_lemma3.4}
\Esub{n\in[N]} \1_{\P_\ell}(n) = 
f(\ell;\mu_N,\sigma_N)\bigg(1 + \Oh\bigg(\frac{1}{\sqrt{\log\log{N}}}\bigg)\bigg).
\end{equation*}
Combined with the Hardy-Ramanujan theorem, this proves that uniformly over all functions $a\colon\N\to\C$ with $\|a\|_\infty\leq 1$,
\begin{align}
\label{eqn_loy_lemma3.4+}
\Esub{n\in[N]} a(\Omega(n))&=\sum_{\ell=1}^{\infty} a(\ell)f(\ell; \mu_N, \sigma_N)+\oh_{N\to\infty}(1).
\end{align}
Juxtaposing \eqref{eqn_Gaussian_density_function_characterization} and \eqref{eqn_loy_lemma3.4+}, it becomes clear that $f(\ell; \mu_N, \sigma_N)$ can be considered the ``probability density function'' of $\Omega(n)$ over the interval $[N]$.
Note that \eqref{eqn_loy_lemma3.4+} is a strengthening of the \Erdos{}-Kac theorem and contains \eqref{eqn_pnt_1} as a special case. It also immediately implies the following reformulation of Theorem~\ref{thm_A}.

\begin{maintheorem}
\label{thm_B}
For any bounded functions $a,b\colon\N\to\C$, we have
\begin{equation*}
    \logE_{n\in[N]} a(\Omega(n))b(\Omega(n+1))
    =
    \sum_{k,\ell=1}^{\infty} a(k)b(\ell)f(k; \mu_N, \sigma_N) f(\ell; \mu_N, \sigma_N)
    + \Oh\bigg(\frac{(\log\log\log N)^2}{\sqrt{\log\log N}}\bigg).
\end{equation*}
where $f(x;\mu,\sigma)$ is the probability density function of the Gaussian distribution defined in \eqref{eqn_Gaussian_prob_dens_fctn}.
Moreover, the error term only depends on $\max\{\|a\|_\infty,\|b\|_\infty\}$.
\end{maintheorem}

As an application of our main theorem, we obtain new results related to the asymptotic behavior of shifted almost primes. In particular, the next result asserts that for any bounded function $a\colon\N\to\C$, and for ``almost all'' $\ell$ in the typical range $T_{A,N}$, the sequence $a(\Omega(n+1))$ with $n$ ranging over $\ell$-almost primes asymptotically behaves the same as the sequence $a(\Omega(n))$ with $n$ ranging over all positive integers.
To state our result in this direction, we define $\overline{\pi}_\ell(N)$ to be the relative density of the set of $\ell$-almost primes within the set of integers up to $N$, that is,
\begin{equation}\label{pi_defn}
    \overline{\pi}_\ell(N)=\frac{|\P_\ell\cap [N]|}{N}.
\end{equation}

\begin{maintheorem}
\label{thm_C}
Let $a\colon\N\to\C$ be a bounded function. Then we have
\begin{equation}\label{eqn1_thm_C}
    \sum_{\ell=1}^\infty \overline{\pi}_\ell(N)\Big|\logE_{n\in[N]\cap\P_\ell} a(\Omega(n+1)) - \Esub{n\in[N]} a(\Omega(n))\Big|
    \ll_{\|a\|_\infty} \Oh\bigg(\frac{(\log\log\log N)^2}{\sqrt{\log\log N}}\bigg).
\end{equation}
In particular, for any fixed $A>1$ and any $\epsilon>0$, for all but at most $\epsilon|T_{A,N}|$ numbers $\ell\in T_{A,N}$, we have 
\begin{equation}\label{eqn2_thm_C}
    \logE_{n\in[N]\cap\P_\ell} a(\Omega(n+1)) = \Esub{n\in[N]} a(\Omega(n)) + \Oh_{\epsilon,\|a\|_\infty}\bigg(\frac{(\log\log\log{N})^2}{\sqrt{\log\log{N}}}\bigg).
\end{equation}
\end{maintheorem}

We remark that if $a(n)=(-1)^n$, then \eqref{eqn1_thm_C} with better bounds than we have is given by \cite[Eq.~(1.16)]{helf_radz}.

We prove that Theorems~\ref{thm_A} and~\ref{thm_C} are equivalent in \cref{Proof_Thm_C}.

We conclude this section with an application of Theorem~\ref{thm_A} to uniform distribution.

\begin{defn}
    Let $(w(n))_{n\in\N}$ be a non-increasing sequence of positive real numbers, and $W(N)=\sum_{n=1}^Nw(n)$ be such that $\lim_{N\to\infty} W(N) =\infty$. A sequence $(\mathbf{x}_n)_{n\in\N}$ in $\R^d$ is {\em uniformly distributed mod $1$ with weights $w(n)$} if 
    $$\lim_{N\to\infty}\frac{1}{W(N)}\sum_{n=1}^N w(n)\1_{\mathbf{I}}(\mathbf{x}_n) = m_{\T^d}(\mathbf{I}),$$
    for any interval $\mathbf{I}\subset\T^d$, where $m_{\T^d}$ denotes the Lebesgue measure on $\T^d$.
\end{defn}

\begin{corollary}\label{ud_corollary}
    Let $\alpha,\beta\in\R$. The sequence $(\alpha\Omega(n),\beta\Omega(n+1))$ is uniformly distributed mod $1$ with weights $\frac{1}{n}$ if and only if both $\alpha,\beta$ are irrational.
\end{corollary}

\begin{proof}
    In view of Weyl's equidistribution criterion, it suffices to show that
    $$\lim_{N\to\infty}\logE_{n\in[N]}e(k\alpha\Omega(n)+\ell\beta\Omega(n+1)) = 0$$
    for all $k,\ell\in\Z$, not both zero, if and only if both $\alpha,\beta$ are irrational. 

    Suppose that at least one of $\alpha,\beta$ is rational, say $\alpha\in\Q$. We write $\alpha = r/q$ for some $r\in\Z$ and $q\in\N$, and we choose $k=q$ and $\ell=0$. Then, for this choice of $k$ and $\ell$, we have
    $$\logE_{n\in[N]} e(k\alpha\Omega(n)+\ell\beta\Omega(n+1)) 
    = \logE_{n\in[N]}e(r\Omega(n)) = 1.$$

    Now suppose that both $\alpha,\beta$ are irrational. In view of Theorem~\ref{thm_A}, it is enough to prove that for any $k\in\Z\setminus\{0\}$, we have 
    $$\lim_{N\to\infty}\Esub{n\in[N]}e(k\alpha\Omega(n))=0.$$
    But this follows from \Halasz's mean value theorem and the fact that the multiplicative function $f(n) := e(k\alpha\Omega(n))$ satisfies $\Dist(f(n),n^{it})=\infty$ for all real $t\neq 0$.\footnote{The definition of the distance $\Dist$ is given in \cref{defn_distance}.} The latter is easy to check, see for instance \cite[Lemmas B.2. \& B.3.]{Charamaras_MVT_24}. This concludes the proof.
\end{proof}

\subsection{Proof strategy for Theorem~\ref{thm_A}}\label{proof strategy}

In broad strokes, we follow the approach introduced by Tao in \cite{Tao16} to prove \cref{thm_log_2_correlation_chowla}, which builds on results of \Matomaki-\Radziwill-Tao in \cite{MRT15}.
The biggest difference in our approach is that the results in \cite{MRT15} cannot be used directly, as they apply exclusively to \emph{strongly aperiodic}\footnote{A multiplicative function $f\colon\N\to\C$ is called \emph{strongly aperiodic} if for any $B>0$, we have $\lim_{N\to \infty}\inf_{|t|\leq BX,q\leq B}\Dist(f,n\mapsto\chi(n) n^{it};N) = \infty$,
where the infimum is taken over all Dirichlet characters $\chi$ of conductor
$q$ with $q\leq B$, and $\Dist$ denotes the \emph{distance} between two multiplicative functions (see \cref{defn_distance}).} multiplicative functions, and the functions we are dealing with here are not even multiplicative.

The starting point in the proof of Theorem~\ref{thm_A} is reducing our problem to the case of multiplicative functions. After replacing the \Cesaro{} averages in the left hand side of \eqref{thm_A_eq} with logarithmic ones, then using Fourier analysis, the \Erdos-Kac theorem and a theorem of Helfgott and \Radziwill{} in \cite{helf_radz}, we can reduce bounding
\begin{equation}\label{eq_pf_strat_1}
    \Big|\logE_{n\in[N]}a(\Omega(n))b(\Omega(n+1)) - \Big(\Esub{n\in[N]}a(\Omega(n))\Big)\Big(\Esub{n\in[N]}b(\Omega(n))\Big)\Big|
\end{equation}
to bounding
\begin{equation}\label{eq_pf_strat_2}
    \sum_{\xi\in I_N}\logE_{n\in[N]}\Big|\logE_{p\in\Pc_N} \f(n+p) - \logE_{m\in[N]}\f(m)\Big|^2,
\end{equation}
where $(\f)_{\xi\in I_N}$ is a sequence of completely multiplicative functions defined in \eqref{defn_f,g}, $I_N$ is a finite integer interval of frequencies $\xi$ defined in \eqref{I_N_def}, the length of which grows with $N$ in appropriate speed, and $\Pc_N$ is a set of primes with size depending on $N$ defined in \eqref{defn_P}. 

Now, if $\f$ was strongly aperiodic for all frequencies $\xi\in I_N$, then one could use the results from \cite{MRT15} (after using a circle method estimate) to bound the expression in question. Note that this would not be enough though, since the sum over $\xi$ could lead the error term to explode, thus, one would need a quantitative lower bound on $\inf_{|t|\leq BX,q\leq B}\Dist(\f,n\mapsto\chi(n) n^{it};N)$, and not just divergence to infinity.
The major issue here is that not all $\f$ are strongly aperiodic. To overcome this, we distinguish between the case when the frequency $\xi$ is small and the case when it is large. For the small frequencies, we can bound the expression in question by establishing quantitative novel bounds for a local \Erdos-Kac theorem. On the other hand, for the large frequencies $\xi$, the function $\f$ is not just strongly aperiodic, but we have convenient lower bound for $\inf_{|t|\leq BX,q\leq B}\Dist(\f,n\mapsto\chi(n) n^{it};N)$. Then, the desired bound in the case of the large frequencies is established by using a circle method estimate and a result from \cite{MRT15}.

In practice, our approach differs slightly from what was described above. For technical reasons, it is more convenient to approximate $\f$ at small values of $\xi$ by another multiplicative function $\g$, defined in \eqref{defn_f,g}, at the very beginning of the proof. Consequently, rather than bounding \eqref{eq_pf_strat_2}, we control \eqref{eq_pf_strat_1} via
\begin{equation}\label{eq_pf_strat_3}
    \sum_{\xi\in I_N}\logE_{n\in[N]}\Big|\logE_{p\in\Pc_N} F_{\xi,N}(n+p) - \logE_{m\in[N]}F_{\xi,N}(m)\Big|^2,
\end{equation}
where $F_{\xi,N}$ equals $\f$ for large $\xi$ and $\g$ for small $\xi$. We carry out this reduction in \cref{reduction_section}, while the contribution from small frequencies is treated in \cref{small_xi_section} and that from large frequencies in \cref{large_xi_section}.

\subsection{Conjectures concerning almost primes}

Establishing, for $h,\ell\in\N$, the asymptotic growth of expressions of the form
\[
\sum_{n=1}^N \1_{\P_\ell}(n)\1_{\P_\ell}(n+h)
\]
is a longstanding and very difficult problem in number theory.
Based on the heuristic that the prime factorization of $n$ and $n+h$ display a high degree of independence, it is reasonable to expect that for all $A>0$, all $k\in\N$ and all $\ell_1,\ldots,\ell_k$ in the typical range $T_{A,N}$, we have
\begin{align}
\label{eqn_almost_primes_independence}
\Esub{n\in[N]} \1_{\P_{\ell_1}}(n)\1_{\P_{\ell_2}}(n+1)\cdots \1_{\P_{\ell_k}}(n+k-1)
&= \overline{\pi}_{\ell_1}(N)\cdots \overline{\pi}_{\ell_k}(N)(1+\oh_{N\to\infty}(1)).
\end{align}
The following conjecture implies that \eqref{eqn_almost_primes_independence} holds for ``almost all'' $\ell_1,\ldots,\ell_k$ in the typical range.

\begin{conjecture}
\label{conj_independence_almost_primes}
For all $k\in\N$, we have
\begin{equation*}
\sum_{\ell_1,\ldots,\ell_k=1}^\infty \Big|\Esub{n\in[N]} \1_{\P_{\ell_1}}(n)\1_{\P_{\ell_2}}(n+1)\cdots \1_{\P_{\ell_k}}(n+k-1) - \overline{\pi}_{\ell_1}(N)\cdots\overline{\pi}_{\ell_k}(N)\Big|
= \oh_{N\to\infty}(1).
\end{equation*}
\end{conjecture}

\cref{conj_independence_almost_primes} implies \cref{conj_functional_chowla}, although the relationship between the two is not immediately apparent. To clarify the connection, we offer below an equivalent formulation of \cref{conj_independence_almost_primes}.

\begin{conjecture}
\label{conj_functional_chowla_str}
For any $k\in\N$ and any bounded $a\colon\N^k\to\C$, we have
$$
\Esub{n\in[N]} a\big(\Omega(n),\Omega(n+1),\ldots,\Omega(n+k-1)\big)
= \Esub{(n_1,\ldots,n_k)\in[N]^k} a\big(\Omega(n_1),\ldots,\Omega(n_k)\big) + \oh_{N\to\infty}(1).
$$
\end{conjecture}

Conjectures \ref{conj_independence_almost_primes} and \ref{conj_functional_chowla_str} are equivalent. The proof of this equivalence is analogous to the proof that Theorems~\ref{thm_A} and~\ref{thm_C} are equivalent given in \cref{Proof_Thm_C}, so it is omitted.

\vspace{2mm} 
\noindent
\textbf{Acknowledgements.}
This work was supported by the Swiss National Science Foundation grant TMSGI2-211214. The first author was supported by the CRM-ISM Postdoctoral Fellowship and the Fondation Courtois - Chaire Courtois en recherche fondamentale - II during the completion of this work.
The authors thank Harald A.~Helfgott and the anonymous referees for valuable comments and suggestions and Biao Wang for pointing out minor mistakes in previous versions of this article.

%SECTION
\section{Preparatory results}\label{prep_lemmas}

In this section we state and prove preliminary results that we will use to prove Theorems~\ref{thm_A} and~\ref{thm_C}.
We start with a simple bound that is used multiple times throughout. For any $C>0$,
\begin{equation*}
    \int_C^\infty e^{-\frac{u^2}{2}} \d u
    \leq \frac{1}{C}\int_C^\infty u e^{-\frac{u^2}{2}} \d u
    = \frac{e^{-\frac{C^2}{2}}}{C}.
\end{equation*}
In particular,
\begin{equation}\label{int_bound}
    \int_{\sqrt{\log\log\log N}}^\infty e^{-\frac{u^2}{2}}\d u \leq \frac{1}{\sqrt{\log\log N}}.
\end{equation}
    
\subsection{Passing from \Cesaro{} to logarithmic averages}

The following lemma gives a well-known identity relating \Cesaro{} averages to logarithmic ones.

\begin{lemma}\label{Cesaro_to_log}
Let $f\colon\N\to\C$ be $1$-bounded. Then for any $N\in\N$ and $\epsilon\in\left(0,\frac{1}{2}\right)$, we have
$$\logE_{n\in[N]} f(n) = \logE_{M\in(N^\epsilon,N)}\Esub{n\in[M]} f(n) + \Oh\bigg(\frac{1}{\log{N}} + \epsilon\bigg).$$
\end{lemma}

\begin{proof}
By partial summation, we have 
\begin{align*}
    \sum_{n\leq N} \frac{f(n)}{n} 
    & = \sum_{M<N}\frac{1}{M(M+1)}\sum_{n\leq M} f(n) + \Esub{n\in[N]} f(n) \\
    & = \sum_{M<N}\frac{1}{M}\Esub{n\in[M]} f(n) + \Esub{n\in[N]} f(n) + \Oh(1),
\end{align*}
hence,
$$\logE_{n\in[N]} f(n) = \logE_{M\in[1,N)}\Esub{n\in[M]} f(n) + \Oh\bigg(\frac{1}{\log{N}}\bigg).$$
Moreover, we have
\begin{align*}
    & \Big|\logE_{M\in(N^\epsilon,N)}\Esub{n\in[M]} f(n) - \logE_{M\in[1,N)}\Esub{n\in[M]} f(n)\Big| \\
    & \hspace*{3cm} \leq \bigg|\frac{1}{\log{N}}\sum_{M\in(N^\epsilon,N)}\frac{1}{M}\Esub{n\in[M]} f(n) - \frac{1}{\log{N}}\sum_{M<N}\frac{1}{M}\Esub{n\in[M]} f(n)\bigg| + \Oh(\epsilon) \\
    & \hspace*{3cm} \leq \frac{1}{\log{N}}\sum_{M\leq N^\epsilon}\frac{1}{M} + \Oh(\epsilon) 
    = \Oh(\epsilon).
\end{align*}
Combining the last two estimates, the lemma follows.
\end{proof}

Recall that given $\ell,N\in\N$ we have defined $\overline{\pi}_\ell(N) = \frac{|\P_\ell\cap [N]|}{N}$. We also define
$$\overline{\pi}^{\log}_\ell(N) = \logE_{n\in[N]}\1_{\P_\ell}(n).$$

We note the following immediate identities that will be invoked several times throughout \cref{Proof_Thm_C}: For any function $a\colon\N\to\C$ and any $N\in\N$,
\begin{equation}\label{eqn_expanding_a(Omega(n))_1}
    a(\Omega(n)) = \sum_{\ell=1}^\infty a(\ell)\1_{\P_\ell}(n),
\end{equation}
and consequently,
\begin{equation}\label{eqn_expanding_a(Omega(n))_2}
    \Esub{n\in[N]} a(\Omega(n)) 
    = \sum_{\ell=1}^\infty a(\ell)\overline{\pi}_\ell(N)
    \quad\text{and}\quad
    \logE_{n\in[N]} a(\Omega(n)) 
    = \sum_{\ell=1}^\infty a(\ell)\overline{\pi}^{\log}_\ell(N).
\end{equation}
   
\begin{lemma}\label{pi_k_estimate}
    For any $A>1$ and $N\in\N$ sufficiently large, we have
    \begin{multline*}
        \overline{\pi}_\ell(N) 
        = f(\ell;\mu_N,\sigma_N)\bigg(1 + \Oh\bigg(\frac{1}{\sqrt{\log\log{N}}}\bigg)\bigg) \\
        = \frac{\exp\Big(-\frac{1}{2}\Big(\frac{\ell-\log\log{N}}{\sqrt{\log\log{N}}}\Big)^2\Big)}{\sqrt{2\pi\log\log{N}}}
    \bigg(1 + \Oh\bigg(\frac{1}{\sqrt{\log\log{N}}}\bigg)\bigg),
    \end{multline*}
    uniformly for all $\ell\in T_{A,N}$, where $T_{A,N}$ is the typical range defined in \eqref{typical_range}, and $f(\ell;\mu_N,\sigma_N)$ is the Gaussian probability density function defined in \eqref{eqn_Gaussian_prob_dens_fctn} with mean $\mu_N$ and standard deviation $\sigma_N$ defined in \eqref{mu_N and sigma_N}.
\end{lemma}

\begin{proof}
    This follows immediately from the estimate
    $$\overline{\pi}_\ell(N) = \frac{\exp\Big(-\frac{1}{2}\Big(\frac{\ell-\log\log{N}}{\sqrt{\log\log{N}}}\Big)^2\Big)}{\sqrt{2\pi\log\log{N}}}
    \bigg(1 + \Oh\bigg(\frac{1}{\sqrt{\log\log{N}}} + \frac{|\ell-\log\log{N}|^3}{(\log\log{N})^2}\bigg)\bigg),$$
    which holds uniformly for all positive integers $\ell\leq \frac{3}{2}\log\log{N}$. For this estimate we refer the reader to \cite[p.~236]{MV_book}.
\end{proof}

\begin{lemma}\label{pk_and_logpk}
    For $N\in\N$ sufficiently large, we have  
    $$\sum_{\ell=1}^\infty \big|\overline{\pi}_\ell(N) - \overline{\pi}^{\log}_\ell(N)\big|
    \ll \frac{1}{\sqrt{\log\log{N}}}.$$
\end{lemma}

\begin{proof}
    By \eqref{int_bound}
    we have
    $$\frac{1}{\sqrt{2\pi}}\int_{-A}^{A} e^{-\frac{u^2}{2}} \d u 
    = 1 + \Oh\bigg(\frac{1}{\sqrt{\log\log{N}}}\bigg).$$
    We define $T_N = (\log\log{N} - A\sqrt{\log\log{N}}, \log\log{N} + A\sqrt{\log\log{N}}]\cap\Z$ and we denote its complement in $\Z$ by $T_N^\text{c}$. Then, by the quantitative \Erdos-Kac theorem \cite[Theorem 2]{renyi-turan}, we have
    $$\Esub{n\in[N]} \1_{T_N^\text{c}}(\Omega(n))
    = 1 - \frac{1}{\sqrt{2\pi}}\int_{-A}^{A} e^{-\frac{u^2}{2}} \d u + \Oh\bigg(\frac{1}{\sqrt{\log\log{N}}}\bigg)
    \ll \frac{1}{\sqrt{\log\log{N}}}.$$
    In view of \cref{Cesaro_to_log}, we also have
    $$\logE_{n\in[N]} \1_{T_N^\text{c}}(\Omega(n))
    \ll \frac{1}{\sqrt{\log\log{N}}}.$$
    By the two last bounds, it follows that
    $$\sum_{\ell\in T_N^\text{c}} \overline{\pi}_\ell(N)
    = \Esub{n\in[N]}\sum_{\ell\in T_N^\text{c}} \1_{\P_\ell}(n)
    = \Esub{n\in[N]} \1_{T_N^\text{c}}(\Omega(n))
    \ll \frac{1}{\sqrt{\log\log{N}}}$$
    and respectively,
    $$\sum_{\ell\in T_N^\text{c}} \overline{\pi}^{\log}_\ell(N)
    \ll \frac{1}{\sqrt{\log\log{N}}}.$$
    Therefore, using triangle inequality, we see that it is enough to prove
    \begin{equation*}
        \sum_{\ell\in T_N}\big|\overline{\pi}_\ell(N) - \overline{\pi}^{\log}_\ell(N)\big|
        \ll \frac{1}{\sqrt{\log\log{N}}}.
    \end{equation*}
    Applying \cref{Cesaro_to_log} with $\epsilon = (\log\log{N})^{-\frac{1}{2}}$ and the triangle inequality, we have
    $$\sum_{\ell\in T_N}\big|\overline{\pi}_\ell(N) - \overline{\pi}^{\log}_\ell(N)\big|
    \leq \logE_{M\in(N^\epsilon,N)}\sum_{\ell\in T_N}\big|\overline{\pi}_\ell(N) - \overline{\pi}_\ell(M)\big| + \Oh\bigg(\frac{1}{\sqrt{\log\log N}}\bigg),$$
    hence, in view of \cref{pi_k_estimate}, it suffices to show that 
    \begin{equation}\label{restriction_TN}
        \logE_{M\in(N^\epsilon,N)}\sum_{\ell\in T_N}|f(\ell;\mu_N,\sigma_N) - f(\ell;\mu_M,\sigma_M)|
        \ll \frac{1}{\sqrt{\log\log{N}}}.
    \end{equation}
    Given $t\in\R$, we define the function $f_t(u) = f(t;u,u)$ and we set $x=\log\log M$, $y=\log\log N$. Then, the left-hand side of \eqref{restriction_TN} becomes 
    $$\logE_{M\in(N^\epsilon,N)}\sum_{\ell\in T_N}|f_\ell(x)-f_\ell(y)|.$$
    By the Fundamental Theorem of Calculus and the triangle inequality, this is
    $$\leq \logE_{M\in(N^\epsilon,N)}\int_x^y\sum_{\ell\in T_N}|f_\ell'(u)| \d u,$$
    where $f_\ell'$ is the derivative of $f_\ell$.
    Now we bound the sum by 
    $$\int_{T_N} |f_t'(u)|\d t + |f'_{t_1}(u)|+|f'_{t_2}(u)|,$$
    where the $t_1,t_2$ are the points in $T_N$ in which the two local maxima of the function $t\mapsto |f'_t(u)|$ are attained, and both maxima are $\Oh(u^{-1})$. The details are left as an exercise to the interested reader. Then, the left hand side of \eqref{restriction_TN} is 
    $$\leq \logE_{M\in(N^\epsilon,N)}\int_x^y\int_{T_N} |f_t'(u)|\d t \d u + \Oh\bigg(\log x - \logE_{M\in(N^\epsilon,N)}\log y\bigg).$$
    To calculate the error term, first we have 
    $$\sum_{N^\epsilon<M<N}\frac{1}{M} = (1-\epsilon)\log N + \Oh(1),$$
    and
    \begin{multline*}
        \sum_{N^\epsilon<M<N}\frac{\log\log\log M}{M} 
        = \int_{N^\epsilon}^N \frac{\log\log\log t}{t} \d t + \Oh(1) 
        = \Big[(\log t)\log\log\log t - \mathrm{li}(\log t)\Big]_{N^\epsilon}^N \\
        = (\log N)\log\log\log N - \mathrm{li}(\log N) - \epsilon(\log N)\log\log\log(N^\epsilon) + \mathrm{li}(\epsilon\log N).
    \end{multline*}
    Using $\mathrm{li}(x)\sim\frac{x}{\log x}$, a simple calculation gives that the error term is 
    $$\frac{\mathrm{li}(\log N)}{(1-\epsilon)\log N} + \frac{\epsilon}{1-\epsilon} \log\left(1 + \frac{\log\epsilon}{\log\log N}\right) - \frac{\mathrm{li}(\epsilon\log N)}{(1-\epsilon)\log N} + \Oh\bigg(\frac{1}{\log N}\bigg) \ll \frac{1}{\log\log N}.$$

    In addition, we can also extend the inner integral to $\R$, thus, in view of \eqref{restriction_TN}, it is enough to show 
    \begin{equation}\label{eq_double_int}
        \logE_{M\in(N^\epsilon,N)}\int_x^y\int_{-\infty}^\infty |f_t'(u)|\d t \d u
        \ll \frac{1}{\sqrt{\log\log N}}.
    \end{equation}
    To show this, we first estimate the inner integral.
    Calculating 
    $$|f_t'(u)| = \frac{|t^2-u^2-u|}{2u^2}f_t(u),$$
    we have
    $$\int_{-\infty}^\infty |f_t'(u)\d t = \frac{1}{2\sqrt{2\pi}}\int_{-\infty}^\infty \frac{|t^2-u^2-u|}{u^2\sqrt{u}}e^{-\frac{(t-u)^2}{2u}}\d t.$$
    Applying the change of variables $t=u+s\sqrt{u}$ and using the triangle inequality, the above is
    $$\leq \frac{1}{\sqrt{2\pi u}}\int_{-\infty}^\infty |s|e^{-\frac{s^2}{2}}\d s + \frac{1}{2u\sqrt{2\pi}}\int_{-\infty}^\infty |s^2-1|e^{-\frac{s^2}{2}}\d s 
    = \frac{\mathbb{E}|Z|}{\sqrt{u}} + \frac{\mathbb{E}|Z^2-1|}{2u},$$
    for $Z\sim \mathcal{N}(0,1)$. Using the well-known estimate $\mathbb{E}|Z|=\sqrt{\frac{2}{\pi}}$ along with the elementary computation $\mathbb{E}|Z^2-1|=\frac{4}{\sqrt{2\pi e}}$, it follows that the inner integral in \eqref{eq_double_int} is $\ll u^{-\frac{1}{2}}$. Then,
    $$\int_x^y \int_{-\infty}^\infty |f_t'(u)|\d t \d u \ll \sqrt{y} - \sqrt{x} = \frac{y-x}{\sqrt{y}+\sqrt{x}}\asymp \sqrt{\log\log N} - \frac{\log\log M}{\sqrt{\log\log N}},$$
    so, the left-hand side of \eqref{eq_double_int} is 
    $$\ll \sqrt{\log\log N} - \frac{1}{\sqrt{\log\log N}}\logE_{M\in(N^\epsilon,N)}\log\log M.$$
    Now we have 
    \begin{align*}
        \sum_{N^\epsilon<M<N}\frac{\log\log M}{M}
        & = \int_{N^\epsilon}^N \frac{\log\log t}{t} \d t + \Oh(1)
        = \Big[(\log t)(\log\log t - 1)\Big]_{N^\epsilon}^N + \Oh(1) \\
        & = (1-\epsilon)(\log N)(\log\log N - 1) + \frac{(\log N)\log\log\log N}{2\sqrt{\log\log N}},
    \end{align*}
    thus,    
    $$\logE_{M\in(N^\epsilon,N)}\log\log M = \log\log N - 1 + \Oh\Bigg(\frac{\log\log\log N}{\sqrt{\log\log N}}\Bigg),$$
    and then \eqref{eq_double_int} follows, concluding the proof.
\end{proof}

\begin{lemma}\label{aOmega_and_logaOmega}
For any bounded function $a\colon \N\to\C$ and $N\in\N$ sufficiently large, we have  
$$
\Big|\Esub{n\in[N]} a(\Omega(n))-\logE_{n\in[N]} a(\Omega(n))\Big| 
\ll_{\|a\|_\infty} \frac{1}{\sqrt{\log\log{N}}}.
$$
\end{lemma}

\begin{proof}
Invoking \eqref{eqn_expanding_a(Omega(n))_2}, we have
$$
\Big|\Esub{n\in[N]} a(\Omega(n))-\logE_{n\in[N]} a(\Omega(n))\Big|
=
\Bigg|\sum_{\ell=1}^\infty a(\ell)\big(
\overline{\pi}_\ell(N) -\overline{\pi}^{\log}_\ell(N)\big)
\Bigg|
\ll
\sum_{\ell=1}^\infty \big|
\overline{\pi}_\ell(N) -\overline{\pi}^{\log}_\ell(N)\big|.
$$
The claim now follows from \cref{pk_and_logpk}.
\end{proof}

\subsection{Independence of periodic functions}

We say that a function $f\colon\N\to\C$ can be written as a \emph{linear combination of $k$ indicator functions of infinite arithmetic progressions} if there exist $v_1,\ldots,v_k\in\{z\in\C: |z|\leq 1\}$, $a_1,\ldots,a_k\in\N$ and $b_1,\ldots,b_k\in\N\cup\{0\}$ such that
$$
f(n)=
v_1 \1_{a_1\N+b_1}(n)+\ldots+v_k \1_{a_k\N+b_k}(n),
%\sum_{i=1}^k c_i \1_{a_i\N+b_i}(n),
\qquad\forall n\in\N.
$$
Note that the modulus of the coefficients $v_1,\ldots,v_k$ in the above description is not allowed to exceed $1$.

\begin{lemma}
\label{lem_independence_periodic_fctns_np_1}
Let $k,\ell,r,s\in\N$. Let $f\colon\N\to\C$ be an $r$-periodic function that can be written as a linear combination of $k$ indicator functions of infinite arithmetic progressions and $g\colon \N\to\C$ an $s$-periodic function that can be written as a linear combination of $\ell$ indicator functions of infinite arithmetic progressions. If $r$ and $s$ are coprime, then for any $N\in\N$,
\begin{equation}
\label{eqn_periodic_fctn_independence_np_1}
\Esub{n\in[N]} f(n)g(n) = \Big(\Esub{n\in[N]} f(n)\Big)\Big(\Esub{n\in[N]} g(n)\Big) + \Oh\bigg(\frac{k\ell}{N}\bigg),
\end{equation}
where the implicit constant in $\Oh(k\ell/N)$ is universal (i.e., does not depend on $k$, $\ell$, $r$, $s$, $f$, $g$, or $N$).
\end{lemma}

\begin{proof}
Let  
$f(n)=v_1 \1_{a_1\N+b_1}(n)+\ldots+v_k \1_{a_k\N+b_k}(n)$ and $g(n)=w_1 \1_{c_1\N+d_1}(n)+\ldots+w_\ell \1_{c_\ell\N+d_\ell}(n)$.
Then
\begin{align}
\label{eqn_periodic_fctn_independence_01}
\Esub{n\in[N]} f(n)g(n)
= \sum_{\substack{1\leq i \leq k\\ 1\leq j \leq \ell}} v_i w_j \Big(\Esub{n\in[N]} \1_{a_i\N+b_i}(n)\1_{c_j\N+d_j}(n)\Big).
\end{align}
Since $a_i$ divides $r$ and $c_j$ divides $s$, we get $\gcd(a_i,c_j)=1$. So 
$$\Esub{n\in[N]} \1_{a_i\N+b_i}(n)\1_{c_j\N+d_j}(n)= \frac{1}{a_i c_j} + \Oh\left(\frac{1}{N}\right).$$ 
On the other hand, we have 
$$\Esub{n\in[N]} \1_{a_i\N+b_i}(n) = \frac{1}{a_i}+\Oh\left(\frac{1}{N}\right)
\quad\text{and}\quad
\Esub{n\in[N]} \1_{c_j\N+d_j}(n) = \frac{1}{c_j} + \Oh\left(\frac{1}{N}\right).$$
This implies that
\begin{align}\label{eqn_periodic_fctn_independence_02}
\Esub{n\in[N]} \1_{a_i\N+b_i}(n)\1_{c_j\N+d_j}(n)
= \Big(\Esub{n\in[N]} \1_{a_i\N+b_i}(n)\Big)\Big(\Esub{n\in[N]} \1_{c_j\N+d_j}(n)\Big) +\Oh\left(\frac{1}{N}\right).
\end{align}
Combining \eqref{eqn_periodic_fctn_independence_01} and \eqref{eqn_periodic_fctn_independence_02} gives
\begin{align*}\label{eqn_periodic_fctn_independence_03}
\Esub{n\in[N]} f(n)g(n)
&=
\sum_{\substack{1\leq i \leq k\\ 1\leq j \leq \ell}} v_i w_j \Big(\Esub{n\in[N]} \1_{a_i\N+b_i}(n)\Big)\Big(\Esub{n\in[N]} \1_{c_j\N+d_j}(n)\Big) +\Oh\left(\frac{k\ell}{N}\right)
\\
&=
\Big(\Esub{n\in[N]} f(n)\Big)\Big(\Esub{n\in[N]} g(n)\Big)+\Oh\left(\frac{k\ell}{N}\right),
\end{align*}
as was to be shown. 
\end{proof}

\subsection{\Turan{}-Kubilius inequality and applications}

The following is a special case of the \Turan{}-Kubilius inequality, which we will use to derive rough upper bounds on the distribution of $\Omega(n)$ outside the typical range. Let $n\mapsto \1_{p\mid n}$ denote the indicator function of the set of numbers divisible by a fixed prime $p$.

\begin{theorem}[\Turan{}-Kubilius inequality; see \cite{elliott_TK}]{}
\label{nam_turan_kubilius}
For any $N\in\N$ and $P\subset\P\cap [1,N]$, we have
\begin{equation}
\label{eqn_turan}
\Esub{n\in[N]}\Bigg|\sum_{p\in P}\1_{p\mid n}-\sum_{p\in P}\frac{1}{p}\Bigg| \leq2\Bigg(\sum_{p\in P}\frac{1}{p}\Bigg)^{\frac{1}{2}}.
\end{equation}
\end{theorem}

\begin{corollary}
\label{cor_TK_1}
For $N\in\N$ and $D\geq 1$, let $$\Gamma_{N,D}=\Bigg\{n\leq N: \sum_{p\leq N}\1_{p\mid n}\geq \log\log N+ D\sqrt{\log\log{N}}\Bigg\}$$ and $$\Delta_{N,D} = \Bigg\{n\leq N\colon \sum_{p\leq N}\1_{p\mid n} \leq \log\log N - D\sqrt{\log\log{N}}\Bigg\}.$$
Then for $N$ sufficiently large, we have
$$\frac{|\Gamma_{N,D}|}{N} =\Oh\bigg(\frac{1}{D}\bigg) \quad\text{and}\quad
\frac{|\Delta_{N,D}|}{N} =\Oh\bigg(\frac{1}{D}\bigg).$$
\end{corollary}

\begin{proof}
We prove the result only for $\Gamma_{N,D}$, since the proof for $\Delta_{N,D}$ is analogous.
By \cref{nam_turan_kubilius} and Markov's inequality, we have for all $t>0$ and all but finitely many $N$ that
\[
\frac{1}{N}\Bigg|\Bigg\{n\leq N: \Bigg|\sum_{p\leq N}\1_{p\mid n}-\sum_{p\leq N}\frac{1}{p}\Bigg| >t \Bigg\}\Bigg|
\leq \frac{2}{t} \Bigg(\sum_{p\leq N}\frac{1}{p}\Bigg)^{\frac{1}{2}}
\leq \frac{4\sqrt{\log\log{N}}}{t}.
\]
It follows that
\[
\frac{1}{N}\Bigg|\Bigg\{n\leq N: \sum_{p\leq N}\1_{p\mid n} \geq \sum_{p\leq N}\frac{1}{p} + t \Bigg\}\Bigg|\leq  \frac{4\sqrt{\log\log{N}}}{t}.
\]
Since, by Mertens' theorem, for all but finitely many $N$ we have
\[
\sum_{p\leq N}\frac{1}{p} + \frac{D}{2}\sqrt{\log\log{N}}\leq \log\log N+ D\sqrt{\log\log{N}},
\]
we can take $t=\frac{D}{2}\sqrt{\log\log{N}}$ and get
\[
\frac{|\Gamma_{N,D}|}{N} =\Oh\bigg(\frac{1}{D}\bigg)
\]
as desired.
\end{proof}

\begin{lemma}
\label{lem_approx_deg_K}
Let $N\in\N$ and $M\in\N$ such that $N^\frac{1}{\log\log{N}}\leq M\leq N^{2\log\log{N}}$. Let also $t,\sigma\in\R$ with $\sigma\geq \frac{1}{2}$, $P\subset \P\cap [1,N]$, and $p,q,K\leq H$ be positive integers. Consider the degree $K$ Taylor approximation of $e(x)$, given by
\[
e_K(x)=\sum_{k=0}^K \frac{(2\pi i x)^k}{k!}.
\]
If $K\geq 200 t (\log\log N)^\sigma$, then for $N$ sufficiently large, we have
\begin{multline}
\label{eqn_approx_deg_K_1}
\Esub{n\in[N]} e\bigg(\frac{t\left(\sum_{r\in P}\1_{r\mid (n+p)}-\1_{r\mid (n+q)}\right)}{\sqrt{\log\log{N}}}\bigg)
\\
=\Esub{n\in[N]}e_K\bigg(\frac{t\left(\sum_{r\in P}\1_{r\mid (n+p)}-\1_{r\mid (n+q)}\right)}{\sqrt{\log\log{N}}}\bigg) + \Oh\bigg(e^{-K}+\frac{1}{(\log\log N)^\sigma}+\frac{H}{N}\bigg)
\end{multline}
and
\begin{multline}
\label{eqn_approx_deg_K_2}
\Esub{\substack{m\in[M] \\ n\in[N]}}e\bigg(\frac{t\left(\sum_{r\in P}\1_{r\mid m}-\1_{r\mid n}\right)}{\sqrt{\log\log{N}}}\bigg)
\\
=\Esub{\substack{m\in[M] \\ n\in[N]}}e_K\bigg(\frac{t\left(\sum_{r\in P}\1_{r\mid m}-\1_{r\mid n}\right)}{\sqrt{\log\log{N}}}\bigg) + \Oh\bigg(e^{-K}+\frac{1}{(\log\log N)^\sigma}\bigg).
\end{multline}
\end{lemma}

\begin{proof}
We only provide a proof of \eqref{eqn_approx_deg_K_1}. The proof of \eqref{eqn_approx_deg_K_2} is analogous. Note that
\begin{align*}
|e(x)-e_K(x)|
=\bigg|\sum_{k=K+1}^\infty \frac{(2\pi i x)^k}{k!}\bigg|
& = \bigg|(2\pi i x)^{K+1}\sum_{k=0}^\infty \frac{(2\pi i x)^{k}}{(k+K+1)!}\bigg|
\\
&\leq \bigg|\frac{(2\pi i x)^{K+1}}{(K+1)!}\bigg| \bigg|\underbrace{\sum_{k=0}^\infty \frac{(2\pi i x)^{k}}{k!}}_{=1}\bigg|.
\end{align*}
Let $D=(\log\log N)^\sigma$, and let $\Gamma_{N,D}$ be as in \cref{cor_TK_1}. Then, by the conclusion of \cref{cor_TK_1}, we have
\begin{align*}
& \Esub{n\in[N]} e\bigg(\frac{t\left(\sum_{r\in P}\1_{r\mid (n+p)}-\1_{r\mid (n+q)}\right)}{\sqrt{\log\log{N}}}\bigg) 
\\
& \hspace*{.5cm} = \Esub{n\in[N]} \1 _{\Gamma_{N,D}}(n+p)\1 _{\Gamma_{N,D}}(n+q) e\bigg(\frac{t\left(\sum_{r\in P}\1_{r\mid (n+p)}-\1_{r\mid (n+q)}\right)}{\sqrt{\log\log{N}}}\bigg) \\
& \hspace*{4cm} + \Oh\bigg(\frac{1}{(\log\log N)^\sigma}+\frac{H}{N}\bigg)
\\
& \hspace*{.5cm} = \Esub{n\in[N]} \1 _{\Gamma_{N,D}}(n+p)\1 _{\Gamma_{N,D}}(n+q) e_K\bigg(\frac{t\left(\sum_{r\in P}\1_{r\mid (n+p)}-\1_{r\mid (n+q)}\right)}{\sqrt{\log\log{N}}}\bigg)
\\
& \hspace*{4cm} + \Oh\bigg(\sup_{n\in \Gamma_{N,D}} \frac{20^{K+1} t^{K+1} \left(\sum_{r\in P}\1_{r\mid n}\right)^{K+1} }{(K+1)!(\log\log N)^{\frac{K+1}{2}}} +\frac{1}{(\log\log N)^\sigma}+\frac{H}{N}\bigg)
\\
& \hspace*{.5cm} = \Esub{n\in[N]}  e_K\bigg(\frac{t\left(\sum_{r\in P}\1_{r\mid (n+p)}-\1_{r\mid (n+q)}\right)}{\sqrt{\log\log{N}}}\bigg)
\\
&\hspace*{4cm} + \Oh\bigg(\sup_{n\in \Gamma_{N,D}} \frac{20^{K+1} t^{K+1} \left(\sum_{r\in P}\1_{r\mid n}\right)^{K+1} }{(K+1)!(\log\log N)^{\frac{K+1}{2}}} +\frac{1}{(\log\log N)^\sigma}+\frac{H}{N}\bigg).
\end{align*}
We now have
\begin{align*}
\sup_{n\in \Gamma_{N,D}} \frac{20^{K+1} t^{K+1} \left(\sum_{r\in P}\1_{r\mid n}\right)^{K+1} }{(K+1)!(\log\log N)^{\frac{K+1}{2}}} 
&\leq 
\frac{40^{K+1} t^{K+1} (\log\log N)^{(K+1)\sigma} }{(K+1)!}.  
\end{align*}
Note, by Stirling's formula, we have $\log(m!)\geq m\log m-2m$ for all but finitely many $m$, and hence for all but finitely many $N$ we have that
\begin{align*}
&\log\bigg(\frac{40^{K+1} t^{K+1} (\log\log N)^{(K+1)\sigma} }{(K+1)!}\bigg)
\\
&\qquad\qquad\leq (K+1)\log(40t)+ (K+1)\sigma(\log\log\log N)- (K+1)\log(K)+2(K+1)
\\
&\qquad\qquad\leq (K+1)\left(\log(40t)+ \sigma (\log\log\log N)- \log(K)+2\right)
\\
&\qquad\qquad= (K+1)\bigg(\log\bigg(\frac{40 t (\log\log N)^\sigma}{K}\bigg)+2\bigg).
\end{align*}
Invoking the assumption $K\geq 200 t (\log\log N)^\sigma$, it now follows that
\[
\log\bigg(\frac{40^{K+1} t^{K+1} (\log\log N)^{\sigma(K+1)} }{(K+1)!}\bigg)\leq -(K+1)
\]
and hence
\begin{align*}
\sup_{n\in \Gamma_{N,D}} \frac{20^{K+1} t^{K+1} \left(\sum_{r\in P}\1_{r\mid n}\right)^{K+1} }{(K+1)!(\log\log N)^{\frac{K+1}{2}}} 
&\leq e^{-(K+1)}.
\end{align*}
The conclusion now follows by combining the above estimates.
\end{proof}

\subsection{Approximations of $\Omega(n)$}

For technical reasons, it will be convenient to replace $\Omega(n)$ with an approximation $\omega_N(n)$ that ignores the contribution of large primes.
We define 
\begin{equation}\label{omega_N_def}
     t_N=N^{\frac{1}{(\log\log N)^8}}
     \quad\text{and}\quad
     \omega_N(n)=\sum_{\substack{p\in\P \\ p\leq t_N}} \1_{p\mid n}.
\end{equation}

The following lemma estimates the error that occurs when replacing $\Omega(n)$ with $\omega_N(n)$ in the expressions appearing in the proof of our main theorem.

\begin{lemma}
\label{lem_approx_Omega}
Let $N,H\geq1$ and $p\leq H$ be positive integers.
For any $t\in\R$, we have
$$\Esub{n\in[N]} e\bigg(\frac{t\Omega(n+p)}{\sqrt{\log\log{N}}}\bigg)
\\
=\Esub{n\in[N]}  e\bigg(\frac{t\omega_N(n+p)}{\sqrt{\log\log{N}}} \bigg)+\Oh\left(\frac{t(\log\log\log N)}{\sqrt{\log\log{N}}}+\frac{tH}{N}\right),$$
and
$$\Esub{n\in[N]} e\bigg(\frac{t\Omega(n)}{\sqrt{\log\log{N}}}\bigg)
= \Esub{n\in[N]} e\bigg(\frac{t\omega_N(n)}{\sqrt{\log\log{N}}}\bigg) + \Oh\left(\frac{t(\log\log\log N)}{\sqrt{\log\log{N}}}\right).$$
\end{lemma}

\begin{proof}
We only prove the first estimate, as the proof of the second is analogous.
Using that $|e(x)-1|^2 = 2-2\cos{x} \leq x^2$ for all $x\in\R$, we get
\begin{align*}
&\Esub{n\in[N]} \left|e\bigg(\frac{t\Omega(n+p)}{\sqrt{\log\log{N}}}\bigg)-e\bigg(\frac{t\omega_N(n+p)}{\sqrt{\log\log{N}}} \bigg)\right|
= \Esub{n\in[N]} \left|e\bigg(\frac{t\left|\Omega(n+p) - \omega_N(n+p)\right|}{\sqrt{\log\log{N}}}\bigg)-1\right|
\\
&\leq \frac{t}{\sqrt{\log\log{N}}}
\Esub{n\in[N]}\left|\Omega(n+p) - \omega_N(n+p)\right|
\ll \frac{t}{\sqrt{\log\log{N}}}\Esub{n\in[N]}\left|\Omega(n)-\omega_N(n)\right| + \frac{tH}{N}.
\end{align*}
Now, we have
\begin{multline*}
\Esub{n\in[N]}\left|\Omega(n)-\omega_N(n)\right|
=\Esub{n\in[N]}\sum_{t_N< p \leq N}\1_{p\mid n}+\Oh(1)
=\sum_{t_N< p \leq N} \frac{1}{N}\left\lfloor \frac{N}{p}\right\rfloor+\Oh(1)
\\
\leq \sum_{t_N< p \leq N} \frac{1}{p} +\Oh(1)
\ll \log\bigg(\frac{\log{N}}{\log{t_N}}\bigg)
= 8\log\log\log{N},
\end{multline*}
concluding the proof.
\end{proof}

%SECTION
\section{Reduction to multiplicative functions}\label{reduction_section}

The next three sections are dedicated to proving Theorem~\ref{thm_A}. In this section, using a key proposition (\cref{Helf-Radz_cor}), which follows from \cite[Corollary 1.4]{helf_radz} and a Fourier analytic argument, we reduce the statement in Theorem~\ref{thm_A} to an assertion concerning averages that involve a sequence of multiplicative functions $(F_{\xi,N})_{\xi\in I_N}$, where $I_N$ is a set of frequencies (see \cref{reducing_thm_A}).

First, note that in light of \cref{aOmega_and_logaOmega}, we may substitute the \Cesaro{} averages in the left-hand side of \eqref{thm_A_eq} by logarithmic ones, hence Theorem~\ref{thm_A} reduces to establishing that for any bounded functions $a,b\colon\N\to\C$, we have
\begin{equation}\label{thm_A_eq_with_log}
    \logE_{n\in[N]} a(\Omega(n))b(\Omega(n+1)) 
    = \Big(\logE_{n\in[N]}a(\Omega(n))\Big)\Big(\logE_{n\in[N]}b(\Omega(n)\Big)
    + \Oh\bigg(\frac{(\log\log\log N)^2}{\sqrt{\log\log N}}\bigg).
\end{equation}
Note that by normalizing, it is enough to prove the above for $1$-bounded functions. Throughout the next three sections we fix two $1$-bounded functions $a,b\colon\N\to\C$ and $N\in\N$ sufficiently large (depending on any absolute constant that will arise).

We define 
\begin{equation}\label{I_N_def}
    A = A(N) := \sqrt{\log\log\log N}
    \quad\text{and}\quad
    I_N := (-A\sqrt{\log\log{N}},A\sqrt{\log\log{N}}]\cap\Z.
\end{equation}
We define the parameter
\begin{equation}\label{defn_B}
    B = B(N) := \frac{2\sqrt{10}}{\pi}A\sqrt{\log\log\log N}
    = \frac{2\sqrt{10}}{\pi}\log\log\log N
    = \frac{2\sqrt{10}}{\pi}A^2
\end{equation}
to distinguish the small and the large frequencies (see \cref{proof strategy}).
For $\xi\in I_N$, we define the multiplicative function $\g$ and the completely multiplicative function $\f$ by
\begin{equation}\label{defn_f,g}
    \g(n) = e\bigg(\frac{\xi\omega_N(n)}{|I_N|}\bigg)
    \quad\text{and}\quad
    \f(n) := e\bigg(\frac{\xi\Omega(n)}{|I_N|}\bigg),
\end{equation}
where we recall the definition of $\omega_N$ in \eqref{omega_N_def}.
We also define, for $\xi\in I_N$,
\begin{equation}\label{defn_F}
    F_{\xi,N} = \g\1_{|\xi|\leq B} + \f\1_{|\xi|>B}.
\end{equation}
For simplicity, we define
$$\Omega_N(n) := \frac{\Omega(n)-\log\log{N}}{\sqrt{\log\log{N}}}.$$
We let
\begin{equation}\label{defn_H}
    H = H(N) := \exp\Big((\log N)^{\frac{1}{\log\log\log N}}\Big) = \exp\bigg(\exp\bigg(\frac{\log\log N}{\log\log\log N}\bigg)\bigg),
\end{equation}
\begin{equation}\label{defn_H0}
    H_0 = H_0(N) := \exp\Big((\log N)^{\frac{3}{4\log\log\log N}}\Big)
    = \exp\bigg(\exp\bigg(\frac{3\log\log N}{4\log\log\log N}\bigg)\bigg)
\end{equation}
and finally,
\begin{equation}\label{defn_P}
    \Pc_N := \P\cap[H_0,H].
\end{equation}
Note that
\begin{equation}\label{eqn_L}
    \mathscr{L} = \mathscr{L}_N := \sum_{p\in\Pc_N}\frac{1}{p} \asymp \frac{\log\log N}{\log\log\log N}.
\end{equation}

We start with a simple tail bound on $\Omega(n)$ that is a consequence of R\'enyi-Tur\'an's quantitative \Erdos-Kac theorem (see \cite[Theorem 2]{renyi-turan}).

\begin{lemma}\label{renyi-turan_cor}
    We have
    $$\frac{|\{n\leq N\colon \Omega(n)>\log\log{N}+A\sqrt{\log\log{N}}\}|}{N}
    \ll \frac{1}{{\sqrt{\log\log{N}}}}$$
    and
    $$\frac{|\{n\leq N\colon \Omega(n)<\log\log{N}-(A-1)\sqrt{\log\log{N}}\}|}{N}
    \ll \frac{1}{{\sqrt{\log\log{N}}}}.$$    
\end{lemma}

\begin{proof}
We prove the first bound, as the proof of the other is analogous.
By \cite[Theorem 2]{renyi-turan}, we have
$$\frac{|\{n\leq N\colon \Omega(n)>\log\log{N}+A\sqrt{\log\log{N}}\}|}{N}
= \frac{1}{\sqrt{2\pi}}\int_A^\infty e^{-\frac{u^2}{2}}\d u + \Oh\bigg(\frac{1}{\sqrt{\log\log{N}}}\bigg).$$
The claim follows by invoking \eqref{int_bound}.
\end{proof}

The next key proposition follows from \cite[Corollary 1.4]{helf_radz}.
\begin{proposition}\label{Helf-Radz_cor}
    We have
    \begin{equation*}
        \logE_{n\in[N]} a(\Omega(n))b(\Omega(n+1))
        = \logE_{p\in\Pc_N}\logE_{n\in[N]} a(\Omega(n)-1)b(\Omega(n+p)-1) + \Oh\Bigg(\sqrt{\frac{\log\log\log N}{\log\log N}}\Bigg).
    \end{equation*}
\end{proposition}

\begin{proof}[Proof of \cref{Helf-Radz_cor}]
    First, we have
    \begin{align*}
        \logE_{n\in[N]} a(\Omega(n)) & b(\Omega(n+1)) \notag \\
        & = \logE_{p\in\Pc_N}\logE_{n\in[N]}  a(\Omega(pn)-1)b(\Omega(pn+p)-1) \notag \\
        & =\logE_{p\in\Pc_N}\logE_{n\in[N]} p\1_{p\mid n}a(\Omega(n)-1)b(\Omega(n+p)-1) + \Oh\bigg(\frac{\log{H}}{\log{N}}\bigg).
    \end{align*}
    Observe that for any $\epsilon>0$, we can restrict the logarithmic averages from over $n$ from $[N]$ to $(N^\epsilon,N]$, gaining an extra $\Oh(\epsilon)$ error term. We choose $\epsilon = (\log\log{N})^{-1}$. Now we split the interval $(N^\epsilon,N]$ into dyadic intervals as follows: we let $J_1 = \left\lfloor\frac{\epsilon(\log{N})}{\log{2}}\right\rfloor+1$, $J_2 = \left\lceil\frac{\log{N}}{\log{2}}\right\rceil$ and then 
    $(N^\epsilon,N] = \bigcup_{J_1\leq j \leq J_2}(2^{j-1},2^j]$. We note that uniformly for all $1$-bounded functions $f$, we have
    \begin{multline*}
        \Esub{j\in[J_1,J_2]}\logE_{n\in(2^{j-1},2^j]} f(n) 
        \asymp
        \frac{\log{2}}{(1-\epsilon)(\log{N})}\sum_{J_1\leq j\leq J_2}\frac{1}{\log{2}}\sum_{2^{j-1}<n\leq 2^j} f(n) \\
        \frac{1}{(1-\epsilon)(\log{N})}\sum_{N^\epsilon<n\leq N} f(n)
        \asymp \logE_{n\in(N^\epsilon,N]} f(n).
    \end{multline*}
    Applying this for $f_p(n)=\logE_{n\in(N^\epsilon,N]}p\1_{p\mid n} a(\Omega(n)-1)b(\Omega(n+p)-1)$, $p\in\Pc_N$, we obtain
    \begin{multline*}
        \logE_{n\in[N]}a(\Omega(n))b(\Omega(n+1)) \\
        = \Esub{j\in[J_1,J_2]}\logE_{p\in\Pc_N}\logE_{n\in(2^{j-1},2^j]}p\1_{p\mid n} a(\Omega(n)-1)b(\Omega(n+p)-1) + \Oh\bigg(\frac{1}{\log\log{N}}\bigg).
    \end{multline*}
   Now we use \cite[Corollary 1.4]{helf_radz} for each $j\in[J_1,J_2]$, to remove the factor $p\1_{p\mid n}$. It is straightforward to check that by the choice of $H_0,H$ and by \eqref{eqn_L}, the assumption of this result is satisfied for all these $j$. Then, we have
    \begin{multline*}
        \logE_{p\in\Pc_N} \logE_{n\in(2^{j-1},2^j]}p\1_{p\mid n}a(\Omega(n)-1)b(\Omega(n+p)-1) \\
        = \logE_{p\in\Pc_N}\logE_{n\in(2^{j-1},2^j]}a(\Omega(n)-1)b(\Omega(n+p)-1) + \Oh\Bigg(\sqrt{\frac{\log\log\log N}{\log\log N}}\Bigg).
    \end{multline*}
    By averaging this over all $j\in[J_1,J_2]$ and combining with the previous, the result follows.
\end{proof}

The next result allows us to reduce the problem to bounding averages involving the multiplicative functions $\g$ and $\f$.

\begin{theorem}\label{reducing_thm_A}
We have
\begin{multline*}
    \Big|\logE_{n\in[N]}a(\Omega(n))b(\Omega(n+1)) - \Big(\logE_{n\in[N]}a(\Omega(n))\Big)\Big(\logE_{n\in[N]}b(\Omega(n))\Big)\Big| \\
    \leq \Bigg(\sum_{\xi\in I_N}\logE_{n\in[N]}\Big|\logE_{p\in\Pc_N} F_{\xi,N}(n+p) - \logE_{m\in[N]}F_{\xi,N}(m)\Big|^2\Bigg)^\frac{1}{2} + \Oh\bigg(\frac{(\log\log\log N)^2}{\sqrt{\log\log N}}\bigg).
\end{multline*}
\end{theorem}

We prove this theorem below in this section.
In view of it, the next two theorems immediately imply \eqref{thm_A_eq_with_log}, thus, Theorem~\ref{thm_A}.

\begin{theorem}\label{small_xi_thm}
We have
\begin{equation*}
    \sum_{|\xi|\leq B}
    \logE_{n\in[N]}\Big|\logE_{p\in\Pc_N}\g(n+p) - \logE_{m\in[N]}\g(m)\Big|^2 
    \ll \frac{(\log\log\log N)^2}{\log\log N}.
\end{equation*}
\end{theorem}

\begin{theorem}\label{large_xi_thm}
We have
\begin{equation*}\label{large_xi_eq}
    \sum_{\substack{\xi\in I_N \\ |\xi|>B}}
    \logE_{n\in[N]}\Big|\logE_{p\in\Pc_N}\f(n+p) - \logE_{m\in[N]}\f(m)\Big|^2 
    \ll \frac{1}{\log\log N}.
\end{equation*}
\end{theorem}

Theorems~\ref{small_xi_thm} and~\ref{large_xi_thm} are proved in Sections~\ref{small_xi_section} and~\ref{large_xi_section} respectively.

\begin{proof}[Proof of \cref{reducing_thm_A}]
We start by defining the following functions:
\begin{itemize}
    \item $\widetilde{a},\widetilde{b}\colon\N\to\C$, given by 
    $$\widetilde{a}(n) 
    = 
    \begin{cases}
        a(n-1), & \text{if}~n\geq 2, \\
        0, & \text{otherwise},
    \end{cases}
    \quad\text{and}\quad
    \widetilde{b}(n) 
    = 
    \begin{cases}
        b(n-1), & \text{if}~n,\geq 2, \\
        0, & \text{otherwise}.
    \end{cases}
    $$
    \item $b_N\colon\N\to\C$ is given by
    $$b_N(m) 
    := 
    \begin{cases}
        b(m+\lfloor\log\log{N}\rfloor), & \text{if}~m\in I_N, \\
        0, & \text{otherwise},
    \end{cases}
    $$
    \item $\widetilde{b}_N\colon\N\to\C$ given by
    $$
    \widetilde{b}_N(m) = 
    \begin{cases}
        b_N(m-1), & \text{if}~m\geq 2 \\
        0, & \text{otherwise}
    \end{cases}
    \,=
    \begin{cases}
        \widetilde{b}(m+\lfloor\log\log{N}\rfloor), & \text{if}~m\in I_N+1 \\
        0, & \text{otherwise.}
    \end{cases}
    $$
\end{itemize}
By \cref{Helf-Radz_cor}, we have
\begin{align}\label{applying_Tao}
    & \Big|\logE_{n\in[N]} a(\Omega(n))b(\Omega(n+1)) - \Big(\logE_{n\in[N]}a(\Omega(n))\Big)\Big(\logE_{n\in[N]}b(\Omega(n))\Big)\Big| \notag \\
    & \hspace*{.8cm} = \Big|\logE_{p\in\Pc_N}\logE_{n\in[N]}\widetilde{a}(\Omega(n))\widetilde{b}(\Omega(n+p)) - \Big(\logE_{n\in[N]} \widetilde{a}(\Omega(n))\Big)\Big(\logE_{n\in[N]}\widetilde{b}(\Omega(n))\Big)\Big| \notag \\
    & \hspace*{10cm} + \Oh\Bigg(\sqrt{\frac{\log\log\log N}{\log\log N}}\Bigg).
\end{align}
We now define the interval 
$$J_N = \left(-A + \frac{1-\{\log\log{N}\}}{\sqrt{\log\log{N}}}, A + \frac{1-\{\log\log{N}\}}{\sqrt{\log\log{N}}}\right].$$
By definition, if $\Omega(n)-\lfloor\log\log{N}\rfloor\in I_N+1$,
which is equivalent to $\Omega_N(n)\in J_N$, then
$\widetilde{b}(\Omega(n)) = \widetilde{b}_N(\Omega(n)-\lfloor\log\log{N}\rfloor)$,
and if $\Omega_N(n)\not\in J_N$, then  $\widetilde{b}_N(\Omega(n)-\lfloor\log\log{N}\rfloor) = 0$. Thus,
\begin{align*}
     &\Big|\logE_{p\in\Pc_N}\logE_{n\in[N]} \widetilde{a}(\Omega(n))\widetilde{b}(\Omega(n+p)) - \Big(\logE_{n\in[N]}\widetilde{a}(\Omega(n))\Big)\Big(\logE_{n\in[N]}\widetilde{b}(\Omega(n))\Big)\Big| \\
     & = \Big|\logE_{p\in\Pc_N}\logE_{n\in[N]}\1_{J_N}(\Omega_N(n+p))\widetilde{a}(\Omega(n))\widetilde{b}_N(\Omega(n+p)-\lfloor\log\log{N}\rfloor) \\
     & \hspace{.6cm} - \Big(\logE_{n\in[N]}\widetilde{a}(\Omega(n))\Big)\Big(\logE_{n\in[N]}\1_{J_N}(\Omega_N(n))\widetilde{b}_N(\Omega(n)-\lfloor\log\log{N}\rfloor)\Big)\Big| 
     + \Oh(R_1 + R_2) \\
     & = \Big|\logE_{p\in\Pc_N}\logE_{n\in[N]}\widetilde{a}(\Omega(n))\widetilde{b}_N(\Omega(n+p)-\lfloor\log\log{N}\rfloor) \\
     & \hspace*{2.65cm} - \Big(\logE_{n\in[N]}\widetilde{a}(\Omega(n))\Big)\Big(\logE_{n\in[N]}\widetilde{b}_N(\Omega(n)-\lfloor\log\log{N}\rfloor)\Big)\Big| 
     + \Oh(R_1 + R_2) \\
     & = \Big|\logE_{n\in[N]}\widetilde{a}(\Omega(n))\Big(\logE_{p\in\Pc_N}\widetilde{b}_N(\Omega(n+p)-\lfloor\log\log{N}\rfloor) \\
     & \hspace*{5.65cm} - \logE_{m\in[N]}\widetilde{b}_N(\Omega(m)-\lfloor\log\log{N}\rfloor)\Big)\Big| + \Oh(R_1+R_2),
\end{align*}
where 
\begin{align*}
    R_1 & = \logE_{p\in\Pc_N}\logE_{n\in[N]} \left(\1_{(A,\infty)}\left(\Omega_N(n+p)\right) + \1_{(-\infty,-A+1)}\left(\Omega_N(n+p)\right)\right)\\
    R_2 & = \logE_{m\in[N]} \left(\1_{(A,\infty)}\left(\Omega_N(m)\right) + \1_{(-\infty,-A+1)}\left(\Omega_N(m)\right)\right).
\end{align*}
By \cref{renyi-turan_cor}, and \cref{Cesaro_to_log} applied for $\epsilon=(\log{N})^{-\frac{1}{2}}$ (so that $N^\epsilon\to\infty$), we can easily see that both $R_1$ and $R_2$ are $\ll(\log\log N)^{-\frac{1}{2}}$.
It follows that
\begin{align}\label{a_to_b_final}
    & \Big|\logE_{p\in\Pc_N}\logE_{n\in[N]}\widetilde{a}(\Omega(n))\widetilde{b}(\Omega(n+p)) - \Big(\logE_{n\in[N]}\widetilde{a}(\Omega(n))\Big)\Big(\logE_{n\in[N]}\widetilde{b}(\Omega(m))\Big)\Big| \notag \\
    & = \Big|\logE_{p\in\Pc_N}\logE_{n\in[N]}\widetilde{b}_N(\Omega(n+p)-\lfloor\log\log{N}\rfloor) - \logE_{n,m\in[N]}\widetilde{b}_N(\Omega(m)-\lfloor\log\log{N}\rfloor)\Big| \notag \\
    & \hspace*{9.7cm}  
    + \Oh\Bigg(\sqrt{\frac{\log\log\log N}{\log\log N}}\Bigg).
\end{align}
Now we apply the Fourier transform to the function $\widetilde{b}_N$ to find
\begin{align*}
    &\widetilde{b}_N(\Omega(m)-\lfloor\log\log{N}\rfloor)
    = b_N(\Omega(m)-\lfloor\log\log{N}\rfloor-1) \\
    & =\sum_{\xi\in I_N}\widehat{b}_N(\xi)e\bigg(\frac{\xi\left(\Omega(m)-\lfloor\log\log{N}\rfloor-1\right)}{|I_N|}\bigg) 
    = \sum_{\xi\in I_N}\widehat{b}_N(\xi)e\bigg(\frac{\xi\left(\lfloor\log\log{N}\rfloor-1\right)}{|I_N|}\bigg)\f(m).
\end{align*}
Substituting the above in the right-hand side of \eqref{a_to_b_final} and rearranging we obtain
\begin{align*}
    & \Big|\logE_{p\in\Pc_N}\logE_{n\in[N]}a(\Omega(n))b(\Omega(n+p)) - \Big(\logE_{n\in[N]}a(\Omega(n))\Big)\Big(\logE_{n\in[N]}b(\Omega(m))\Big)\Big| \\
    & = \Bigg|\sum_{\xi\in I_N} e\bigg(\frac{\xi\left(\lfloor\log\log{N}\rfloor-1\right)}{|I_N|}\bigg)\widehat{b}_N(\xi) \\
    & \hspace*{1.4cm} \logE_{n\in[N]}\widetilde{a}(\Omega(n))
    \Big(\logE_{p\in\Pc_N}\f(n+p) - \logE_{m\in[N]}\f(m)\Big)\Bigg| + \Oh\Bigg(\sqrt{\frac{\log\log\log N}{\log\log N}}\Bigg).
\end{align*}
Now we split the sum in two pieces; one over the small frequencies $|\xi|\leq B$ and one over the large frequencies $|\xi|>B$.
Using \cref{lem_approx_Omega} to approximate $\f$ by $\g$ in the sum over the small frequencies, we may replace $\f$ by $\g$ at a cost of an error 

$$\Oh\Bigg(\frac{\log\log\log N}{\sqrt{\log\log N}}\sum_{|\xi|\leq B}|\xi\widehat{b}_N(\xi)|\Bigg).$$

Using Cauchy-Schwarz inequality, and then Parseval's identity that gives $\sum_{\xi\in I_N}|\widehat{b}_N(\xi)|^2
= \Esub{n\in I_N}|b_N(n)|^2
\leq \|b\|_\infty
\leq 1$, we can bound the sum in the error by $B\asymp\log\log\log N$.
It follows
\begin{align*}
    & \Big|\logE_{p\in\Pc_N}\logE_{n\in[N]}a(\Omega(n))b(\Omega(n+p)) - \Big(\logE_{n\in[N]}a(\Omega(n))\Big)\Big(\logE_{n\in[N]}b(\Omega(m))\Big)\Big| \notag \\
    & = \Bigg|\sum_{\xi\in I_N} e\bigg(\frac{\xi\left(\lfloor\log\log{N}\rfloor-1\right)}{|I_N|}\bigg)\widehat{b}_N(\xi)\logE_{n\in[N]}\widetilde{a}(\Omega(n)) \notag \\
    & \hspace*{3cm} 
    \Big(\logE_{p\in\Pc_N}F_{\xi,N}(n+p) - \logE_{m\in[N]}F_{\xi,N}(m)\Big)\Bigg| + \Oh\Bigg(\frac{(\log\log\log N)^2}{\sqrt{\log\log{N}}}\Bigg),
\end{align*}
where $F_{\xi,N}$ is defined in \eqref{defn_F}.
Finally, we bound the main term of the last equation. Using Cauchy-Schwarz inequality, we see that it is
$$\leq \Bigg(\sum_{\xi\in I_N}|\widehat{b}_N(\xi)|^2\Bigg)^\frac{1}{2}
\Bigg(\sum_{\xi\in I_N}\Big|\logE_{n\in[N]}\widetilde{a}(\Omega(n))
\Big(\logE_{p\in\Pc_N}F_{\xi,N}(n+p) - \logE_{m\in[N]}F_{\xi,N}(m)\Big)\Big|^2\Bigg)^\frac{1}{2}.$$
Now using Parseval's identity as before we bound the first term of the latter by $1$, while for the second term we apply Cauchy-Schwarz inequality again, using also that $\|a\|_\infty\leq1$, to obtain the overall bound
$$\Bigg(\sum_{\xi\in I_N}\logE_{n\in[N]}
\Big|\logE_{p\in\Pc_N}F_{\xi,N}(n+p) - \logE_{m\in[N]}F_{\xi,N}(m)\Big|^2\Bigg)^\frac{1}{2}.$$
This concludes the proof of the theorem.
\end{proof}

To establish Theorem~\ref{thm_A} it remains to prove Theorems~\ref{small_xi_thm} and~\ref{large_xi_thm}. This is carried out in the following two sections respectively.

%SECTION
\section{Small frequencies}\label{small_xi_section}

The goal of this section is to prove \cref{small_xi_thm}, which we first reduce to the following result.

\begin{theorem}\label{small_xi_Cesaro_thm}
Let $1\leq H_1\leq N$, $200\leq A_0 \leq \sqrt{\log\log{N}}$, $1\leq B_0\ll A_0\sqrt{\log\log N}$, and $\Pc\subset\P\cap[1,H_1]$. For $\xi\in[-B_0,B_0]$, we define the multiplicative function
$$\widetilde{g}(n) = \widetilde{g}_{\xi,N,A_0}(n) = e\bigg(\frac{\xi\omega_N(n)}{A_0\sqrt{\log\log{N}}}\bigg).$$
Then, for any $N^{\frac{1}{\log\log{N}}}\leq M\leq N^{2\log\log{N}}$, we have
\begin{equation}\label{eqn_small_xi_Cesaro_thm}
    \sum_{|\xi|\leq B_0} \Esub{n\in[N]}\Big|\logE_{p\in\Pc}\widetilde{g}(n+p)-\Esub{m\in[M]}\widetilde{g}(m)\Big|^2 
    \ll \frac{B_0 \exp\Big(\frac{B_0^2 \log\log H_1}{A_0^2 \log\log N}\Big)}{\log\log H_1} + \frac{B_0H_1}{N}.
\end{equation}
\end{theorem}

\begin{proof}[Proof of \cref{small_xi_Cesaro_thm} $\Longrightarrow$ \cref{small_xi_thm}]
Using \cref{Cesaro_to_log} for $\epsilon=(\log\log N)^{-1}$, Cauchy-Schwarz inequality, and noting that for $M\in(N^\epsilon,N)$ it holds $\log\log{M}\asymp\log\log{N}$, we have 
\begin{align}\label{M,M'_eq1}
    & \sum_{|\xi|\leq B}\logE_{n\in[N]}\Big|\logE_{p\in\Pc_N} \g(n+p) - \logE_{n\in[N]} \g(m)\Big|^2 \notag \\
    & = \sum_{|\xi|\leq B}\logE_{M\in(N^\epsilon,N)}\Esub{n\in[M]}\bigg|\logE_{p\in\Pc_N} \g(n+p) - \logE_{M'\in(N^\epsilon,N)}\Esub{m\in[M']}\g(m)\bigg|^2 \notag \\
    & \hspace*{.5cm} + \Oh\bigg(\frac{B}{\log\log{N}}\bigg) \notag \\
    & \leq \sup_{N^\epsilon<M,M'<N}\sum_{|\xi|\leq B} \Esub{n\in[M]}\Big|\logE_{p\in\Pc_N} \g(n+p) - \Esub{m\in[M']} \g(m) \Big|^2 + \Oh\bigg(\frac{B}{\log\log{N}}\bigg).
\end{align}
For $M,M'\in(N^\epsilon,N)$, we want to apply \cref{small_xi_Cesaro_thm} (with $M$ in the place of $N$ and $M'$ in the place of $M$) to bound the last sum. We choose $A_0=A$, $B_0 = B$, $H_1 = H(M)$ and $\Pc=\Pc_N$. It is straightforward to check that this choice of parameters satisfy the assumptions of \cref{small_xi_Cesaro_thm}. By the choice of $A_0$, for each $\xi\in[-B,B]$, we can approximate $\widetilde{g}=\widetilde{g}_{\xi,N,A}$ by $\g$ at a cost of a negligible error, so that \eqref{eqn_small_xi_Cesaro_thm} remains true if we replace $\widetilde{g}$ with $\g$.
For for any such $M,M'$, by \cref{small_xi_Cesaro_thm} and since $\log\log{H}\asymp\log\log{N}$, we have
\begin{equation*}
    \sum_{|\xi|\leq B} \Esub{n\in[M]}  \Big|\logE_{p\in\Pc_N} \g(n+p) - \Esub{m\in[M']} \g(m) \Big|^2
    \ll \frac{B \exp\Big(\frac{B^2 \log\log H}{A^2 \log\log N}\Big)}{\log\log H}
    \ll \frac{(\log\log\log N)^2}{\log\log N}.
\end{equation*}
Combining this with \eqref{M,M'_eq1} the result follows.
\end{proof}

\begin{proof}[Proof of \cref{small_xi_Cesaro_thm}]
We let $H_1,A_0,B_0,\Pc,M$ and $\widetilde{g}$ be as in the assumptions of the theorem and we define 
$$S = \sum_{|\xi|\leq B_0} \Esub{n\in[N]}\Big|\logE_{p\in\Pc}\widetilde{g}(n+p)-\Esub{m\in[M]}\widetilde{g}(m)\Big|^2.$$

By expanding the square, we have
\begin{align*}
S & = \sum_{|\xi|\leq B_0}\Esub{n\in[N]}\bigg(\logE_{p,q\in\Pc}\widetilde{
g}(n+p)\overline{\widetilde{g}(n+q)} 
\\
& \hspace*{3cm} - 2\Re\Big(\Big(\logE_{p\in\Pc}\overline{\widetilde{g}(n+p)}\Big)\Big(\Esub{m\in[M]}\widetilde{g}(m)\Big)\Big) + \Esub{m,\ell\in[M]}\widetilde{g}(m)\overline{\widetilde{g}(\ell)}\bigg) 
\\
& = \sum_{|\xi|\leq B}\Bigg(
\logE_{\substack{p,q\in\Pc \\ p\neq q}} \Esub{n\in[N]}\widetilde{g}(n+p)\overline{\widetilde{g}(n+q)}
\\
& \hspace*{3cm} - 2\Re\bigg(\Esub{\substack{m\in[M] \\ n\in[N]}}\overline{\widetilde{g}(n)}\widetilde{g}(m)\bigg) + \Esub{m,\ell\in[M]}\widetilde{g}(m)\overline{\widetilde{g}(\ell)}\Bigg) 
+ \Oh\bigg(\frac{B_0}{\log H_1}+\frac{B_0H_1}{N}\bigg) 
\\
& \leq \sum_{|\xi|\leq B_0}\Bigg|
\logE_{\substack{p,q\in\Pc \\ p\neq q}} \Esub{n\in[N]}\widetilde{g}(n+p)\overline{\widetilde{g}(n+q)}
\\
& \hspace*{3cm} - 2\bigg(\Esub{\substack{m\in[M] \\ n\in[N]}}\overline{\widetilde{g}(n)}\widetilde{g}(m)\bigg) + \Esub{m,\ell\in[M]}\widetilde{g}(m)\overline{\widetilde{g}(\ell)}\Bigg| 
+ \Oh\bigg(\frac{B_0}{\log H_1}+\frac{B_0H_1}{N}\bigg) 
\end{align*}
Let
$K = \lfloor(\log\log{N})^5\rfloor$ and
$e_K(x)=\sum_{k=0}^K \frac{(2\pi i)^k x^k}{k!}$
be the degree $K$ approximation of $e(x)$.
Since $B_0\leq \sqrt{\log\log{N}}$ and $A_0\geq 200$, we have $K\geq \frac{200 B_0 (\log\log{N})^4}{A_0}$. So using the definitions of $\widetilde{g}$ and $\omega_N$, \cref{lem_approx_deg_K}, and that $\log\log{M}\asymp\log\log{N}$, then the main term above is equal to 
\begin{multline*}
    \sum_{|\xi|\leq B_0}\Bigg|
    \logE_{\substack{p,q\in\Pc \\ p\neq q}} \Esub{n\in[N]}e_K\bigg(\frac{\xi \big(\sum_{r\leq t_N}\1_{r\mid n+p}-1_{r\mid n+q}\big)}{A_0\sqrt{\log\log{N}}}\bigg) - 2\Esub{\substack{m\in[M] \\ n\in[N]}}e_K\bigg(\frac{\xi\big(\sum_{r\leq t_N}\1_{r\mid m}-\1_{r\mid n}\big)}{A_0\sqrt{\log\log{N}}}\bigg) \\
    + \Esub{m,\ell\in[M]} e_K\bigg(\frac{\xi\big(\sum_{r\leq t_N}\1_{r\mid m}-\1_{r\mid \ell}\big)}{A_0\sqrt{\log\log{N}}}\bigg)\Bigg| + \Oh\bigg(B_0 e^{-K}+\frac{B_0}{(\log\log N)^4}+\frac{B_0H_1}{N}\bigg),
\end{multline*}
where the sums over $r$ are sums over primes, and throughout this proof the variable $r$ will denote a prime. By definition of $e_K$, the above is 
\begin{align*}
    & \ll \sum_{\xi\leq B_0}\sum_{k=1}^K \frac{\xi^k}{A_0^k k! (\log\log N)^{\frac{k}{2}}}\Bigg|\logE_{\substack{p,q\in\Pc \\ p\neq q}} \Esub{n\in[N]} \bigg(\sum_{r\leq t_N}\1_{r\mid n+p}-\1_{r\mid n+q}\bigg)^k \\
    & \hspace{4.5cm} - 2\Esub{\substack{m\in[M] \\ n\in[N]}}
    \bigg(\sum_{r\leq t_N} \1_{r\mid m}-\1_{r\mid n}\bigg)^k + \Esub{m,\ell\in[M]}\bigg(\sum_{r\leq t_N}\1_{r\mid m}-\1_{r\mid \ell}\bigg)^k\Bigg| \\
    & \ll \sum_{k=1}^K \frac{B_0^{k+1}}{A_0^k(k+1)!(\log\log{N})^\frac{k}{2}}\Bigg|\logE_{\substack{p,q\in\Pc \\ p\neq q}} \Esub{n\in[N]} \bigg(\sum_{r\leq t_N}\1_{r\mid n+p}-\1_{r\mid n+q}\bigg)^k \\
    & \hspace{4.5cm} - 2\Esub{\substack{m\in[M] \\ n\in[N]}}
    \bigg(\sum_{r\leq t_N} \1_{r\mid m}-\1_{r\mid n}\bigg)^k + \Esub{m,\ell\in[M]}\bigg(\sum_{r\leq t_N}\1_{r\mid m}-\1_{r\mid \ell}\bigg)^k\Bigg|. 
\end{align*}
We define $S_1(k)$ to the above absolute value in the last equation, and then we have
\begin{equation}\label{eqn_S1}
    S \ll \sum_{k=1}^K \frac{B_0^{k+1}}{A_0^k(k+1)!(\log\log{N})^\frac{k}{2}}S_1(k) + \Oh\bigg(\frac{B_0}{\log{H_1}} + \frac{B_0H_1}{N} + \frac{B_0}{(\log\log{N})^4}\bigg).
\end{equation}

Given $u_1,\ldots,u_j\geq 0$ with $u_1+\ldots+u_j=k$, let
$$\binom{k}{u_1,\ldots,u_j}$$
denote the multinomial coefficient, which expresses the number of unique permutations of a word of length $k$ using $j$ distinct letters such that the $i$-th letter occurs $u_i$-any times. 
By expanding, we have
\begin{align*}
\bigg(\sum_{r\leq t_N}\1_{r\mid n+p} & - \1_{r\mid n+q}\bigg)^k
= \sum_{r_1,\ldots,r_k\leq t_N} \prod_{i=1}^k \left(\1_{r_i\mid n+p} - \1_{r_i\mid n+q}\right) \\
& = \sum_{j=1}^k \sum_{\substack{u_1,\ldots,u_j\geq 1 \\ u_1+\ldots+u_j=k}}\binom{k}{u_1,\ldots,u_j}\sum_{r_1<\ldots<r_j\leq t_N}\prod_{i=1}^j \left(\1_{r_i\mid n+p} - \1_{r_i\mid n+q}\right)^{u_i} \\
& = 
\sum_{j=1}^k \sum_{\substack{u_1,\ldots,u_j\geq 1 \\ u_1+\ldots+u_j=k}}\binom{k}{u_1,\ldots,u_j}\sum_{r_1<\ldots<r_j\leq t_N}\prod_{i=1}^j \left(\1_{r_i\mid n+p} - \1_{r_i\mid n+q}\right)^{[u_i]},
\end{align*}
where $[u]_2=1$ if $u$ is odd and $[u]_2=2$ if $u$ is even.
In view of \cref{lem_independence_periodic_fctns_np_1}, for any primes $r_1<\ldots<r_j$, letting $\epsilon = (\log\log{N})^{-1}$, we have
$$\Esub{n\in[N]}\prod_{i=1}^{j}\left(\1_{r_i\mid n+p}-\1_{r_i\mid n+q}\right)^{[u_i]_2} 
= \prod_{i=1}^{j}\Esub{n\in[N]}\left(\1_{r_i\mid n+p}-\1_{r_i\mid n+q}\right)^{[u_i]_2}
+ \Oh\bigg(\frac{2^{j^2}}{N}\bigg),$$
$$\Esub{\substack{m\in[M] \\ n\in[N]}} \prod_{i=1}^{j} \left(\1_{r_i\mid m}-\1_{r_i\mid n}\right)^{[u_i]_2}
= \prod_{i=1}^{j}\Esub{\substack{m\in[M] \\ n\in[N]}} \left(\1_{r_i\mid m}-\1_{r_i\mid n}\right)^{[u_i]_2}
+ \Oh\bigg(\frac{2^{j^2}}{N^\epsilon}\bigg)$$
and 
$$\Esub{m,\ell\in[M]} \prod_{i=1}^{j} \left(\1_{r_i\mid m}-\1_{r_i\mid \ell}\right)^{[u_i]_2}
= \prod_{i=1}^{j}\Esub{m,\ell\in[M]}\left(\1_{r_i\mid m}-\1_{r_i\mid \ell}\right)^{[u_i]_2}
+ \Oh\bigg(\frac{2^{j^2}}{N^\epsilon}\bigg).$$
It follows that
\begin{align*}
&\Esub{n\in[N]}\bigg(\sum_{r\leq t_N}\1_{r\mid n+p}-\1_{r\mid n+q}\bigg)^k
\\
& = \sum_{j=1}^k \sum_{\substack{u_1,\ldots,u_j\geq 1 \\ u_1+\ldots+u_j=k}} \binom{k}{u_1,\ldots,u_{t_N}} \sum_{r_1<\ldots<r_j\leq t_N }\Bigg(\prod_{i=1}^j \Esub{n\in[N]}  
\left(\1_{r_i\mid n+p}-\1_{r_i\mid n+q}\right)^{[u_i]_2} + \Oh\bigg(\frac{2^{j^2}}{N^\epsilon}\bigg)\Bigg) \\
& = \sum_{j=1}^k \sum_{\substack{u_1,\ldots,u_{j}\geq 1 \\ u_1+\ldots+u_j=k}} \binom{k}{u_1,\ldots,u_j}\underbrace{\sum_{r_1<\ldots<r_j\leq t_N}\prod_{i=1}^{j} \Esub{n\in[N]}  
\left(\1_{r_i\mid n+p}-\1_{r_i\mid n+q}\right)^{[u_i]_2}}_{:=\Phi_1(u_1,\ldots,u_j;p,q)} + \Oh\bigg(\frac{2^{k^2}t_N^k}{N^\epsilon}\bigg),
\end{align*}
where for the calculation of the error term in the last step above we used that
\begin{align*}
    \sum_{j=1}^k \sum_{\substack{u_1,\ldots,u_j\geq 1 \\ u_1+\ldots+u_j=k}} \binom{k}{u_1,\ldots,u_{t_N}} \sum_{r_1<\ldots<r_j\leq t_N} 2^{j^2}
    & \leq 2^{k^2} \sum_{j=1}^k \binom{t_N}{j} \sum_{\substack{u_1,\ldots,u_j\geq 1 \\ u_1+\ldots+u_j=k}} \binom{k}{u_1,\ldots,u_{t_N}} \\
    & = 2^{k^2}\sum_{\substack{u_1,\ldots,u_{t_N}\geq 0 \\ u_1+\dots +u_{t_N}=k}}\binom{k}{u_1,\ldots,u_{t_N}}
    = 2^{k^2}t_N^k,
\end{align*}
where the last equality follows from the multinomial theorem.
Similar calculations show 
\begin{multline*}
\Esub{\substack{m\in[M] \\ n\in[N]}}\bigg(\sum_{r\leq t_N}\1_{r\mid m}-\1_{r\mid n}\bigg)^k
\\
=
\sum_{j=1}^k\sum_{\substack{u_1,\ldots,u_j\geq 1 \\ u_1+\ldots+u_j=k}} \binom{k}{u_1,\ldots,u_j} \underbrace{\sum_{r_1<\ldots<r_j\leq t_N}\prod_{i=1}^j\Esub{\substack{m\in[M] \\ n\in[N]}} 
\left(\1_{r_i\mid m}-\1_{r_i\mid n}\right)^{[u_i]_2}}_{:=\Phi_2(u_1,\ldots,u_j)} + \Oh\bigg(\frac{2^{k^2}t_N^k}{N^\epsilon}\bigg)
\end{multline*}
and
\begin{multline*}
\Esub{m,\ell\in[M]}\bigg(\sum_{r\leq t_N}\1_{r\mid m}-\1_{r\mid \ell}\bigg)^k
\\
=
\sum_{j=1}^k\sum_{\substack{u_1,\ldots,u_j\geq 1 \\ u_1+\ldots+u_j=k}} \binom{k}{u_1,\ldots,u_j} \underbrace{\sum_{r_1<\ldots<r_j\leq t_N}\prod_{i=1}^j \Esub{m,\ell\in[M]}  
\left(\1_{r_i\mid m}-\1_{r_i\mid \ell}\right)^{[u_i]_2}}_{:=\Phi_3(u_1,\ldots,u_j)} + \Oh\bigg(\frac{2^{k^2}t_N^k}{N^\epsilon}\bigg).
\end{multline*}
Thus, by the definition of $S_1(k)$ and the triangle inequality, we have
\begin{align}\label{eqn_S1_2}
    S_1(k) \leq \sum_{j=1}^k\sum_{\substack{u_1,\ldots,u_j\geq 1 \\ u_1+\ldots+u_j=k}} \binom{k}{u_1,\ldots,u_j}\bigg| \logE_{\substack{p,q\in\Pc \\ p\neq q}} \Phi_1(&u_1,\ldots,u_j;p,q) - 2\Phi_2(u_1,\ldots,u_j) \notag \\
    & + \Phi_3(u_1,\ldots,u_j)\bigg| + \Oh\bigg(\frac{2^{k^2}t_N^k}{N^\epsilon}\bigg).
\end{align}

For $u_1,\ldots,u_j\geq 1$ with $u_1+\ldots+u_j = k$,
we distinguish two cases:
\begin{itemize}
    \item If $u_i$ is odd for some $1\leq i\leq j$, then
    $$\bigg|\Esub{n\in[N]}\1_{r_i\mid n+p}-\1_{r_i\mid n+q}\bigg| < \frac{2}{N}, \quad
    \bigg|\Esub{\substack{m\in[M] \\ n\in[N]}}\1_{r_i\mid m}-\1_{r_i\mid n}\bigg|<\frac{1}{N^\epsilon}$$
    and 
    $$\bigg|\Esub{m,\ell\in[M]}\1_{r_i\mid m}-\1_{r_i\mid \ell}\bigg|<\frac{1}{N^\epsilon}.$$
    It follows that
    \begin{multline}\label{Phi_bound1}
        \bigg|\logE_{\substack{p,q\in\Pc \\ p\neq q}}\Phi_1(u_1,\ldots,u_j;p,q)-2\Phi_2(u_1,\ldots,u_j)+\Phi_3(u_1,\ldots,u_j)\bigg| \\
        \ll \frac{1}{N^\epsilon}\sum_{r_1<\ldots<r_j\leq t_N} 1
        = \frac{1}{N^\epsilon}\binom{\pi(t_N)}{j}
        < \frac{1}{N^\epsilon}\binom{t_N}{j}. 
    \end{multline}
    The contribution of these terms in \eqref{eqn_S1_2} is 
    \begin{equation}\label{eqn_1contr_S1}
        \ll \frac{1}{M}\sum_{j=1}^k \binom{t_N}{j} \sum_{\substack{u_1,\ldots,u_j\geq 1 \\ u_1+\ldots+u_j=k}} \binom{k}{u_1,\ldots,u_{t_N}} 
        = \frac{1}{N^\epsilon} \sum_{\substack{u_1,\ldots,u_{t_N}\geq 0 \\ u_1+\ldots +u_{t_N}=k}}\binom{k}{u_1,\ldots,u_{t_N}} 
        = \frac{t_N^k}{N^\epsilon},
    \end{equation}
where the last equality follows again from the multinomial theorem.
\item On the other hand, if $u_i$ is even for all $1\leq i\leq j$, then
    \begin{align*}
        \Esub{n\in[N]}\big(\1_{r_i\mid n+p}-\1_{r_i\mid n+q}\big)^{[u_i]_2}
        & = \frac{2}{r_i}\big(1-\1_{r_i\mid p-q}\big)+\Oh\bigg(\frac{1}{N}\bigg), \\
        \Esub{\substack{m\in[M] \\ n\in[N]}}\big(\1_{r_i\mid m} - \1_{r_i\mid n}\big)^{[u_i]_2} 
        & = \frac{2}{r_i} + \Oh\bigg(\frac{1}{N^\epsilon}\bigg), \\
        \Esub{m,\ell\in[M]}\big(\1_{r_i\mid m}-\1_{r_i\mid \ell}\big)^{[u_i]_2}
        & = \frac{2}{r_i}+\Oh\bigg(\frac{1}{N^\epsilon}\bigg).
    \end{align*}
    Taking product over $1\leq i \leq j$ yields
    \begin{align*}
        \prod_{i=1}^j \Esub{n\in[N]}\big(\1_{r_i\mid n+p}-\1_{r_i\mid n+q}\big)^{[u_i]_2}
        & = \prod_{i=1}^j \frac{2}{r_i}\big(1-\1_{r_i\mid p-q}\big) + \Oh\Bigg(\sum_{i=1}^j \frac{1}{N^i}\binom{j}{j-i}\Bigg) \\
        & = \prod_{i=1}^j \frac{2}{r_i}\big(1-\1_{r_i\mid p-q}\big) + \Oh\bigg(\frac{j^2}{N}\bigg)
    \end{align*}
    and similar estimates hold for the other two terms. Thus, we obtain
    \begin{align*}
        & \bigg|\logE_{\substack{p,q\in\Pc \\ p\neq q}} 
        \Phi_1(u_1,\ldots,u_j;p,q)-2\Phi_2(u_1,\ldots,u_j) + \Phi_3(u_1,\ldots,u_j)\bigg| \\
        & \hspace*{1cm} = \logE_{\substack{p,q\in\Pc \\ p\neq q}}\sum_{r_1<\ldots<r_j\leq t_N}\prod_{i=1}^j\frac{2}{r_i}\1_{r_i\mid p-q} + \Oh\Bigg(\frac{j^2\binom{t_N}{j}}{N^\epsilon}\Bigg) \\
        & \hspace*{1cm} = 2^j\sum_{r_1<\ldots<r_j\leq H_1}\frac{1}{r_1\cdots r_j}\logE_{\substack{p,q\in\Pc \\ p\neq q}}\1_{r_1\cdots r_j\mid p-q} + \Oh\Bigg(\frac{j^2\binom{t_N}{j}}{N^\epsilon}\Bigg) \\
        & \hspace*{1cm} \ll 2^j\sum_{r_1<\ldots<r_j\leq H_1} \frac{1}{(r_1\cdots r_j)^2} + \frac{2^j}{\log\log{H_1}}\sum_{r_1<\ldots<r_j\leq H_1}\frac{1}{r_1\cdots r_j} + \frac{j^2\binom{t_N}{j}}{N^\epsilon} \\
        & \hspace*{1cm} \ll \frac{2^j}{j! \log\log{H_1}} + \frac{j^2\binom{t_N}{j}}{N^\epsilon}.
    \end{align*}
    Since $u_1,\ldots,u_j$ are even, notice that $j\leq \frac{k}{2}$. Then the contribution of these terms in \eqref{eqn_S1_2} is 
    \begin{align*}
        & \ll \sum_{j\leq\frac{k}{2}} \Bigg(\frac{2^j}{j!\log\log H_1} + \frac{j^2\binom{t_N}{j}}{N^\epsilon}\Bigg) \sum_{\substack{u_1,\ldots,u_j\geq 1 \notag \\ u_1+\ldots+u_j = \left\lfloor\frac{k}{2}\right\rfloor}}\binom{\left\lfloor\frac{k}{2}\right\rfloor}{u_1,\ldots,u_j} \\
        & \leq (\log\log{H_1})^{\frac{k}{2}-1}\sum_{j\leq\frac{k}{2}}\frac{1}{j!}\sum_{\substack{u_1,\ldots,u_j\geq 1 \\ u_1+\ldots+u_j = \left\lfloor\frac{k}{2}\right\rfloor}}\binom{\left\lfloor\frac{k}{2}\right\rfloor}{u_1,\ldots,u_j} + \frac{k^2}{N^\epsilon} \sum_{\substack{u_1,\ldots,u_{t_N}\geq 0 \\ u_1+\ldots+u_{t_N} = k}}\binom{k}{u_1,\ldots,u_{t_N}}.
    \end{align*}
    The sum in the second term is equal to $t_N^k$ by the multinomial theorem.
    For the first term, if $\genfrac{\{}{\}}{0pt}{}{n}{j}$ denotes the Stirling number of the second kind and $B_n$ the $n$-th Bell number, satisfying $\sum_{j=1}^n \genfrac{\{}{\}}{0pt}{}{n}{j}=B_n$, then 
    $$\sum_{j\leq\frac{k}{2}}\frac{1}{j!}\sum_{\substack{u_1,\ldots,u_j\geq 1 \\ u_1+\ldots+u_j = \left\lfloor\frac{k}{2}\right\rfloor}}\binom{\left\lfloor\frac{k}{2}\right\rfloor}{u_1,\ldots,u_j} = \sum_{j\leq\frac{k}{2}}\genfrac{\{}{\}}{0pt}{}{\left\lfloor\frac{k}{2}\right\rfloor}{j} = B_{\left\lfloor\frac{k}{2}\right\rfloor}.$$
    Using the bound $B_n \ll \big(\frac{n}{\log n}\big)^n$ (see \cite[p.~81]{De_Bruijn_81}), it follows that the total contribution of the terms corresponding to even $u_1,\ldots,u_j$ in \eqref{eqn_S1_2} is
    \begin{align}\label{eqn_2contr_S1}
         \ll \bigg(\frac{k}{2\log k}\bigg)^\frac{k}{2}(\log\log H_1)^{\frac{k}{2}-1} + \frac{k^2 t_N^k}{N^\epsilon}.
    \end{align}
\end{itemize}
Substituting the bounds \eqref{eqn_1contr_S1} and \eqref{eqn_2contr_S1} in \eqref{eqn_S1_2}, we obtain 
$$S_1(k) \ll \bigg(\frac{k}{\log k}\bigg)^{\frac{k}{2}}(\log\log{H_1})^{\frac{k}{2}-1} + \frac{2^{k^2}t_N^k}{N^\epsilon}.$$
Then, by \eqref{eqn_S1}, we get
\begin{multline*}
    S \leq
    \frac{B_0}{\log\log{H_1}} \sum_{k=3}^K \frac{1}{(k+1)!}\Bigg(\frac{B_0}{A_0}\cdot\sqrt{\frac{\log\log H_1}{\log\log N}}\cdot\sqrt{\frac{k}{\log k}}\Bigg)^k \\
    + \frac{2^{K^2} t_N^K B_0}{N^\epsilon} \sum_{k=1}^K \frac{1}{(k+1)!}\bigg(\frac{B_0}{A_0\sqrt{\log\log{N}}}\bigg)^k + \frac{B_0}{\log H_1} + \frac{B_0H_1}{N} + \frac{B_0}{(\log\log N)^4}.
\end{multline*}
It essentially remains to bound the two sums in the last equation. Since $B_0\leq A_0\sqrt{\log\log N}$, then the second sum is trivially bounded by $\Oh(1)$. To bound the first sum we set $x=\frac{B_0^2 \log\log H_1}{A_0^2 \log\log N}$ and by Stirling's formula we bound the sum by 
$$\sum_{k=3}^K\bigg(\frac{x}{k\log k}\bigg)^{k/2} = \sum_{k=3}^K \exp\bigg(\frac{k}{2}(\log x- \log k - \log\log k)\bigg).$$
We define the function $f(t) = \exp\big(\frac{t}{2}(\log x- \log t - \log\log t)\big)$, for $t\geq3$. This function has a maximum at $t_\ast$, where  
$\log t_\ast + \log\log t_\ast + \frac{1}{\log t_\ast} - \log x + 1 =0,$
and we also let $t_0$ given by $t_0\log t_0 = 2x$. It is not hard to check that $t_\ast<t_0$, $t_\ast\sim\frac{x}{e\log x}$ and $t_0\sim\frac{2x}{\log x}$.
Then the contribution of the terms corresponding to $k>t_0$ is trivially bounded by the tail of the a convergent geometric series, so it is negligible, while the sum up to $t_0$ is bounded by $t_0f(t_\ast)$, and it quite straightforward to check that this is 
$$\ll \exp\bigg(\frac{x}{2e\log{x}}(1+\oh(1))\bigg)\frac{x}{\log{x}}
\ll\exp\bigg(\frac{x}{\log{x}}\bigg)\leq \exp\bigg(\frac{B_0^2 \log\log H_1}{A_0^2 \log\log N}\bigg).$$ 
This bound is not as sharp as possible, but it is more than enough for our purpose. 
It follows that
\begin{multline*}
    S \ll \frac{B_0 \exp\Big(\frac{B_0^2 \log\log H_1}{A_0^2 \log\log N}\Big)}{\log\log H_1} + \frac{2^{K^2} t_N^K B_0}{N^\epsilon} + \frac{B_0}{\log H_1} + \frac{B_0H_1}{N} + \frac{B_0}{(\log \log N)^4} \\
    \ll \frac{B_0 \exp\Big(\frac{B_0^2 \log\log H_1}{A_0^2 \log\log N}\Big)}{\log\log H_1} + \frac{B_0H_1}{N},
\end{multline*}
finishing the proof.
\end{proof}

%SECTION
\section{Large frequencies}\label{large_xi_section}

The goal of this subsection is to prove \cref{large_xi_thm}.
Expanding the square in 
$$\sum_{\substack{\xi\in I_N \\ |\xi|>B}}
\logE_{n\in[N]}\Big|\logE_{p\in\Pc_N}\f(n+p) - \logE_{m\in[N]} \f(m)\Big|^2$$
we see that \cref{large_xi_thm} follows from the next two propositions.

\begin{proposition}\label{large_xi_thm1}
We have
    $$\sum_{\substack{\xi\in I_N \\ |\xi|>B}}\Big|\Esub{n\in[N]} \f(n)\Big| 
    \ll \frac{1}{(\log\log N)^{50}}.$$
\end{proposition}

\begin{proposition}\label{large_xi_thm2}
We have
    $$\sum_{\substack{\xi\in I_N \\ |\xi|>B}}
    \Esub{n\in[N]}\Big|\logE_{p\in\Pc_N}\f(n+p)\Big|^2 
    \ll \frac{1}{\log\log{N}}.$$
\end{proposition}

We use these propositions for logarithmic averages over $n$, but we can prove them for \Cesaro{} averages, which, in view of \cref{Cesaro_to_log}, makes the results stronger.

To prove \cref{large_xi_thm1} we use \Halasz's mean value theorem, while the proof of \cref{large_xi_thm2} requires using a deeper result of \Matomaki-\Radziwill-Tao \cite{MRT15}. For both cases we need to introduce the notion of distance between two multiplicative functions.

\begin{defn}\label{defn_distance}
    Let $f,g\colon\N\to\C$ be two bounded multiplicative functions. The \textit{distance (up to $N$)} between $f$ and $g$ is defined as $$\Dist(f,g;N)=\Bigg(\sum_{p\leq N}\frac{1-\Re(f(p)\overline{g(p)})}{p}\Bigg)^{\frac{1}{2}}.$$ 
\end{defn}

\subsection{Distance from twisted characters}

To prove Propositions \ref{large_xi_thm1} and \ref{large_xi_thm2}, we need to have appropriate lower bounds on the distance of the functions $\f$ from twisted characters, i.e., functions of the form $\chi(n)n^{it}$, where $\chi$ is a Dirichlet character and $t\in\R$.

\begin{lemma}\label{dist}
Let $\xi\in I_N$ and $t\in\R$.
\begin{itemize}
    \item[(i)] If $|t|\ll1$, then 
    \begin{equation}\label{dist1.eq1}
        \Dist(\f,n\mapsto n^{it};N)^2 
        = \bigg(1-\cos\bigg(\frac{2\pi\xi}{|I_N|}\bigg)\bigg)\log\log{N} 
        + \cos\bigg(\frac{2\pi\xi}{|I_N|}\bigg)\log(1+|t|\log{N}) + \Oh(1).
    \end{equation}
    In particular, for $N$ sufficiently large and for $\xi\in I_N$, we have 
    \begin{equation}\label{dist1.eq2}
        \inf_{|t|\ll 1}\Dist(\f,n\mapsto n^{it};N)^2\geq 
        \begin{cases}
            \frac{\pi^2\xi^2}{2A^2} + \Oh(1), & |\xi|\leq\sqrt{|I_N|}, \\
            \frac{\pi^2 \log\log N}{A} + \Oh(1), & \sqrt{|I_N|}<|\xi|\leq \frac{|I_N|}{4}, \\
            \log\log N + \Oh(1), & \text{otherwise}.
        \end{cases}
    \end{equation}
    \item[(ii)] If $2\leq |t|\leq N^B$, for $B>0$, then for any $\epsilon>0$, we have
    \begin{equation}\label{dist2.eq1}
        \Dist(\f,n\mapsto n^{it}; N)^2\geq \bigg(\frac{1}{3}-\epsilon\bigg)\log\log{N}+\Oh_{B,\epsilon}(1).
    \end{equation}
\end{itemize}
\end{lemma}

\begin{lemma}\label{dist with chi}
Let $\xi\in I_N$, $H'\geq 2$, $\chi$ be a non-principal Dirichlet character modulo $q$, and $t\in\R$.
\begin{itemize}
    \item[(i)] If $q\leq(\log{H'})^5$ and $|t|\ll1$, then 
    $$\Dist(\f,n\mapsto \chi(n)n^{it};N)^2 \gg \log\log{N}.$$
    \item[(ii)] If $q\leq (\log{N})^\frac{1}{125}$ and $10\leq |t|\leq N^B$, for $B>0$, then for any $\epsilon>0$, we have
    $$\Dist(\f,n\mapsto \chi(n)n^{it}; N)^2 \geq \bigg(\frac{1}{3}-\epsilon\bigg)\log\log{N} + \Oh_{B,\epsilon}(1).$$
\end{itemize}
\end{lemma}

\begin{proof}[Proof of \cref{dist}]
(i) Let $s=1+\frac{1}{\log{N}}+it$, $|t|\ll1$, and for $\Re(s)>1$, define $$F(s)=\sum_{n\geq1}\frac{\f(n)}{n^s}.$$ Then we have
\begin{equation}\label{L-f}
    \log{F(s)}=\sum_{p}\frac{\f(p)}{p^s}+\Oh(1) = e\bigg(\frac{\xi}{|I_N|}\bigg)\log{\zeta(s)} + \Oh(1).
\end{equation}
Now, using the following classic estimate for the Riemann zeta function
$$\zeta(s)=\frac{s}{s-1}+O\bigg(\frac{|s|}{\sigma}\bigg)=\frac{1}{s-1}+\Oh(1),
\quad \sigma=\Re(s)>0,$$
we get that
\begin{align}\label{lz}
    \log{\zeta(s)} & =-\log\bigg(\frac{1}{\log{N}}+it\bigg)+\Oh(1) \\
    & = -\log\bigg|\frac{1}{\log{N}}+it\bigg| -  i\arg\bigg(\frac{1}{\log{N}}+it\bigg)+\Oh(1) \notag \\
    & = \log\log{N} - \log(1+|t|\log{N}) + \Oh(1). \notag
\end{align}
Combining \eqref{L-f} and \eqref{lz}, we obtain that
\begin{equation}\label{ReF}
    \Re(\log{F(s)})
    = \cos\bigg(\frac{2\pi\xi}{|I_N|}\bigg)(\log\log{N}-\log(1+|t|\log{N})) + \Oh(1)
\end{equation}
Now, using the asymptotic (see \cite[Eq. (2.2)]{GS_large_char_sums})
$$\exp(-\Dist(\f,n\mapsto n^{it};N)^2)=\frac{|F(s)|}{\log{N}} + \Oh(1),$$ 
we see that
$$\Dist(\f,n\mapsto n^{it};N)^2=\log\log{N}-\log|F(s)|=\log\log{N}-\Re(\log{F(s)})$$ and then, by substituting \eqref{ReF} to the right-hand side of the above we obtain \eqref{dist1.eq1}.
\par
We now prove the second part. We distinguish two cases depending on the size of $|\xi|$:
\begin{itemize}
    \item If $|\xi|\leq|I_N|/4$, then $\cos(2\pi\xi/|I_N|)\geq0$, hence \eqref{dist1.eq1} gives
    $$\Dist(\f,n\mapsto n^{it};N)^2 \geq \bigg(1-\cos\bigg(\frac{2\pi\xi}{|I_N|}\bigg)\bigg)\log\log{N} + O(1).$$
    $1-\cos(2\pi\xi/|I_N|)$ is increasing as a function of $|\xi|$ and it is also positive. Moreover, by the Laurent expansion of $1-\cos(1/x)$ at $x=\infty$, we have
    \begin{equation}\label{Laurent_exp}
        1-\cos\bigg(\frac{1}{x}\bigg) = \frac{1}{2x^2} + \Oh\bigg(\frac{1}{x^4}\bigg).
    \end{equation}
    \begin{itemize}
        \item  If $|\xi|\leq \sqrt{|I_N|}$, then \eqref{Laurent_exp} gives that
        for $N$ sufficiently large, 
        $$1-\cos\bigg(\frac{2\pi\xi}{|I_N|}\bigg) 
        = \frac{2\pi^2\xi^2}{|I_N|^2} + \Oh\bigg(\frac{1}{|I_N|^2}\bigg),$$
        hence 
        $$\Dist(\f,n\mapsto n^{it};N)^2 
        \geq \frac{\pi^2\xi^2}{2A^2} + \Oh(1).$$
        \item If $\sqrt{|I_N|}<|\xi|\leq|I_N|/4$, then by monotonicity and \eqref{Laurent_exp}, we have
        $$1-\cos\bigg(\frac{2\pi\xi}{|I_N|}\bigg) > 1 - \cos\bigg(\frac{2\pi}{\sqrt{|I_N|}}\bigg) = \frac{2\pi^2}{|I_N|} + \Oh\bigg(\frac{1}{|I_N|^2}\bigg),$$
        thus,
        $$\Dist(\f,n\mapsto n^{it};N)^2 \geq \frac{\pi^2(\log\log{N})}{A} + \Oh(1).$$
    \end{itemize}
    \item If $|I_N|/4\leq |\xi| \leq |I_N|/2$, then $\cos(2\pi\xi/|I_N|)\leq 0$. Rewriting \eqref{dist1.eq1}, we have
    $$\Dist(\f,n\mapsto n^{it};N)^2
    = \log\log{N} + \cos\bigg(\frac{2\pi\xi}{|I_N|}\bigg)\log\bigg(\frac{1+|t|\log{N}}{\log{N}}\bigg) + \Oh(1)$$
    and now we distinguish cases on the size of $|t|$:
    If $|t|\leq 1/2$, then the second term above is non-negative for $N$ sufficiently large. On the other hand if $1/2<|t|\ll 1$, then the second term is $\Oh(1)$. Hence, in both cases we have
    $$\Dist(\f,n\mapsto n^{it};N)^2 \geq \log\log{N} + \Oh(1).$$
\end{itemize}
Putting everything together, we get that \eqref{dist1.eq2} holds for $N$ sufficiently large.

(ii) We have that
\begin{align*}
    \Dist(\f,n\mapsto n^{it};N)^2 &= \sum_{p\leq N}\frac{1-\Re\big(e(\xi/|I_N|)p^{-it}\big)}{p} \\
    &\geq \sum_{\exp((\log{N})^{\frac{2}{3}+\epsilon})\leq p\leq N} \frac{1-\Re\big(e(\xi/|I_N|))p^{-it}\big)}{p} \\
    & \geq \sum_{\exp((\log{N})^{\frac{2}{3}+\epsilon})\leq p\leq N}\frac{1}{p}-\Bigg|\sum_{\exp((\log{N})^{\frac{2}{3}+\epsilon})\leq p\leq N}\frac{e(\xi/|I_N|)p^{-it}}{p}\Bigg| \\ 
    &= \bigg(\frac{1}{3}-\epsilon\bigg)\log\log{N}-\Bigg|\sum_{\exp((\log{N})^{\frac{2}{3}+\epsilon})\leq p\leq N}\frac{1}{p^{1+it}}\Bigg|+\Oh(1),
\end{align*}
and so, it suffices to show that $$\Bigg|\sum_{\exp((\log{N})^{\frac{2}{3}+\epsilon})\leq p\leq N}\frac{1}{p^{1+it}}\Bigg|\ll 1.$$ 
The left-hand side above is equal to 
\begin{equation}\label{VK_simple}
    \bigg|\log\bigg(\zeta\bigg(1+\frac{1}{\log{N}}+it\bigg)\bigg) - \log\bigg(\zeta\bigg(1+\frac{1}{(\log{N})^{\frac{2}{3}+\epsilon}}+it\bigg)\bigg)\bigg|+\Oh(1).
\end{equation}
By the Vinogradov-Korobov zero-free region \cite{Korobov_zero-free,Vinogradov_zero-free}, there exists some absolute constant $C>0$, such that $\zeta(\sigma+it)\neq0$ in the region
$$\sigma\geq 1-\frac{C}{(\log|t|)^{\frac{2}{3}}(\log\log|t|)^{\frac{1}{3}}}, \quad |t|\geq 10,$$
and then
it follows using standard analytic methods
(as in \cite{richert}) that \eqref{VK_simple} is $\Oh(1)$.
\end{proof}

\begin{proof}[Proof of \cref{dist with chi}]
(i) For $\Re(s)>1$, we define 
$$F(s,\chi)=\sum_{n\geq1}\frac{\f(n)\overline{\chi(n)}}{n^s}.$$ 
Let $s=1+\frac{1}{\log{N}}+it$, $|t|\ll1$. Then we have
\begin{equation}\label{L-f with chi}
    \log{F(s,\chi)} = \sum_{p}\frac{\f(p)\overline{\chi(p)}}{p^s} + \Oh(1) 
    = e\bigg(\frac{\xi}{|I_N|}\bigg)\log{L(s,\overline{\chi})} + \Oh(1).
\end{equation}
Using the classical asymptotic identity \cite[Eq.~(2.2)]{GS_large_char_sums} as we did in the proof of \cref{dist} (i), we have that
$$\exp(-\Dist(\f,n\mapsto \chi(n)n^{it};N)^2)=\frac{|F(s,\chi)|}{\log{N}} + \Oh(1),$$ 
and combing it with \eqref{L-f with chi}, we obtain that
\begin{align}\label{dist expansion with chi}
    \Dist(\f,&n\mapsto\chi(n)n^{it};N)^2
    = \log\log{N} - \log|F(s,\chi)| + \Oh(1) \notag \\
    & = \log\log{N} - \Re(\log{F(s,\chi)}) + \Oh(1) \notag \\
    & = \log\log{N} - \cos\bigg(\frac{2\pi\xi}{|I_N|}\bigg)\Re(\log{L(s,\overline{\chi})})
    -\sin\bigg(\frac{2\pi\xi}{|I_N|}\bigg)\Im(\log{L(s,\overline{\chi}})) + \Oh(1) \notag \\
    & = \log\log{N} - \cos\bigg(\frac{2\pi\xi}{|I_N|}\bigg)\log|L(s,\overline{\chi})|
    - \sin\bigg(\frac{2\pi\xi}{|I_N|}\bigg)\arg(L(s,\overline{\chi})) + \Oh(1) \notag \\
    & \geq \log\log{N} - \log|L(s,\overline{\chi})| + \Oh(1).
\end{align}
Now we use the following classical bound for Dirichlet L-functions
$$|L(s,\overline{\chi})|\ll \log(q(|t|+2)),$$
which holds uniformly for all $q>1$, all non-principal Dirichlet characters $\chi$
of modulus $q$, all $\sigma\geq1$ and all $t\in\R$ (see for example \cite[Theorem 8.18]{Tenenbaum_book}).
By the assumptions on $q$ and $t$, the above yields 
$$\log|L(s,\overline{\chi})|\ll \log\log\log{H'}.$$
Plugging this in \eqref{dist expansion with chi} yields that
$$\Dist(\f,n\mapsto \chi(n)n^{it};N)^2
\gg \log\bigg(\frac{\log{N}}{\log\log{H'}}\bigg) \gg \log\log{N}.$$
This concludes the proof.

(ii) We omit this proof since it is identical to that of \cref{dist} (ii), with the only difference being that we use the Vinogradov-Korobov zero-free region for Dirichlet $L$-functions (see \cite[p.~176]{Montogomery_ten}), instead of the classical one for the Riemann zeta function. Note that this zero-free region applies here since the fact that $q\leq (\log N)^{1/125}$ ensures the absence of exceptional zeros (see \cite[p.~6]{MRT15}).
\end{proof}

\subsection{Proof of \cref{large_xi_thm1}}

We will use a version of \Halasz's mean value theorem \cite{Hal} due to Tenenbaum \cite{Tenenbaum_book}.
Given a $X,T\geq1$, and $f$ a $1$-bounded multiplicative function, we define 
$$M_0(f;X,T) := \inf_{|t|\leq T} \Dist(f, n\mapsto n^{it};X)^2.$$

\begin{theorem}[{\cite[Theorem 4.5]{Tenenbaum_book}}]\label{Halasz}
For $X\geq3$, $1\leq T\leq (\log X)^3$, and $f$ $1$-bounded multiplicative function, we have 
$$\Esub{n\in[X]} f(n) \ll (1+M_0(f;X,T))e^{-M_0(f,X,T)} + \frac{1}{T} + \frac{1}{\sqrt{\log X}}.$$
\end{theorem}

\begin{proof}[Proof of \cref{large_xi_thm1}]
By \cref{Halasz} for $X=N$ and $T=\log N$, we have
$$\sum_{\substack{\xi\in I_N \\ |\xi|>B}}\Big|\Esub{n\in[N]} \f(n)\Big|
\ll
\sum_{\substack{\xi\in I_N \\ |\xi|>B}}
\exp\bigg(-\frac{M_0(\f;N,\log N)}{2}\bigg) + \frac{|I_N|}{\sqrt{\log N}}.
$$
By \cref{dist}, we have for $N$ sufficiently large,
$$M_0(\f;N,\log N)
\geq \min\bigg\{\frac{\pi^2\xi^2}{2A^2}, \frac{\pi^2(\log\log{N})}{A}\bigg\} + \Oh(1).$$
Then we have
\begin{align*}
    \sum_{\substack{\xi\in I_N \\ |\xi|>B}}\Big|\Esub{n\in[N]} \f(n)\Big|
    \ll \max\Bigg\{\sum_{\xi>B}e^{-\frac{\pi^2\xi^2}{A^2}}, \frac{|I_N|}{(\log{N})^\frac{\pi^2}{16A}}\Bigg\}.
\end{align*}
For $N$ sufficiently large,
\begin{multline*}
    \sum_{\xi>B}e^{-\frac{\pi^2\xi^2}{4A^2}} 
    \leq \int_{B-1}^\infty e^{-\frac{\pi^2 x^2}{4A^2}} \d x 
    = \frac{\sqrt{2}A}{\pi} \int_{{\frac{\pi(B-1)}{\sqrt{2}A}}}^\infty e^{-\frac{y^2}{2}} \d y 
    \leq \frac{2A^2}{\pi^2(B-1)} \int_{{\frac{\pi(B-1)}{\sqrt{2}A}}}^\infty ye^{-\frac{y^2}{2}} \d y \\
    = \frac{2A^2}{\pi^2(B-1)}\exp\bigg(-\frac{\pi^2(B-1)^2}{4A^2}\bigg)
    \ll \exp\bigg(-\frac{\pi^2B^2}{4A^2}\bigg) = (\log\log N)^{-50}.
\end{multline*}
We can readily check that the second term in the maximum is smaller than the above bound, and that the contribution of $\frac{|I_N|}{\sqrt{\log N}}$ in the first equation is negligible, hence the result follows.
\end{proof}

\subsection{Proof of \cref{large_xi_thm2}}

We start with the following circle method estimate.

\begin{lemma}\label{circle_method}
Let $H_1\in\N$, $\Pc\subset[1,H_1]$ be a set of primes, and $(x_j)_{1\leq j\leq H_1}$, $(y_j)_{1-H_1\leq j\leq H_1-1}$ with $x_j,y_j\in\C$. 
For any $\epsilon>0$, we consider the set of major arcs 
$$\mathfrak{M}_{\epsilon,\Pc} = \Big\{\alpha\in[0,1]\colon \Big|\logE_{p\in\Pc}e(p\alpha)\Big|>\epsilon\Big\},$$
and then we have
$$\logE_{p\in\Pc}\Esub{j\in[H_1]} x_j\overline{y_{j-p}} 
\ll
\int_{\mathfrak{M}_{\epsilon,\Pc}}\Big|\Esub{j\in[H_1]} x_je(j\alpha)\Big|\d\alpha + \epsilon.$$
\end{lemma}

\begin{proof}
This is \cite[Lemma 8.1]{helf_radz} (see also \cite[Lemma 5.11]{Helf_Ubis}).
\end{proof}

We also need the following lemma concerning the Lebesgue measure of the major arcs in \cref{circle_method}.

\begin{lemma}\label{major_arcs_bound}
Let $H_1\geq 4$, $\Pc\subset[\frac{H_1}{2},H_1]$ be a set of primes. Then for any $\epsilon>0$, we have
$$|\mathfrak{M}_{\epsilon,\Pc}| \ll \frac{1}{\epsilon^4 H_1},$$
where $\mathfrak{M}_{\epsilon,\Pc}$ is as in \cref{circle_method}.
\end{lemma}

\begin{proof}
    This is \cite[Lemma 8.2]{helf_radz} (see also \cite[Lemma 5.12]{Helf_Ubis}).
\end{proof}

The next lemma essentially allows us to reduce the proof of \cref{large_xi_thm2} to bounding averages of $\f$ twisted by a linear phase in short intervals.

\begin{lemma}\label{bound_dyadic_set}
    Let $f,g\colon\Z\to\C$ be $1$-bounded functions, $H_1\geq 4$, $\Pc\subset[\frac{H_1}{2},H_1]$ be a set of primes. Then for any $\epsilon>0$, we have
    \begin{multline*}
        \logE_{n\in[N]}\logE_{p\in\Pc} f(n)\overline{g(n-p)} \\
        \ll \epsilon^{-4}\bigg(\sup_{\alpha\in[0,1]}\logE_{n\in[N]}\Big|\Esub{j\in[H_1]}f(n+j)e(j\alpha)\Big| + \frac{\log{H_1}}{\log{N}}\bigg) + \epsilon + \frac{\log{H_1}}{\log{N}}.
    \end{multline*}
\end{lemma}

\begin{proof}
    First we have that
    $$\logE_{n\in[N]}\logE_{p\in\Pc} f(n)\overline{g(n-p)}
    = 
    \logE_{n\in[N]}\logE_{p\in\Pc}\Esub{j\in[H_1]}f(n+j)\overline{g(n+j-p)} + \Oh\bigg(\frac{\log{H_1}}{\log{N}}\bigg).$$
    We suppose that $f,g$ are both supported on $[N]$, and this will add an extra error term of order $\frac{\log{H_1}}{\log{N}}$. Then applying \cref{circle_method} for $f(n+j) = x_j$ and $g(n+j) = y_j$, $n\in[N]$, and then \cref{major_arcs_bound}, we obtain        
    \begin{multline*}
        \logE_{n\in[N]}\logE_{p\in\Pc}\Esub{j\in[H_1]}f(n+j)\overline{g(n+j-p)} \\
        \ll \epsilon^{-4}\sup_{\alpha\in[0,1]}\logE_{n\in[N]}\Big|\Esub{j\in[H_1]}f(n+j)e(j\alpha)\Big| + \epsilon + \frac{\log{H_1}}{\log{N}}.
    \end{multline*}
    Dropping the assumption that $f$ is supported on $[N]$ adds an extra error term of order $\epsilon^{-4}\cdot\frac{\log{H_1}}{\log{N}}$.
\end{proof}

Finally, to bound the averages of $\f$ twisted by a linear phase in short intervals, we shall employ a powerful result of \Matomaki-\Radziwill-Tao \cite{MRT15}.
Given $X,H'\geq1$, we let 
$$Q_{X,H'} := \min((\log{X})^\frac{1}{125},(\log{H'})^5)$$
and given a $1$-bounded multiplicative function $f$, we define
$$M(f;X,Q_{X,H'}) := \inf_{\substack{\chi\:(q) \\ q\leq Q_{X,H'}}}\inf_{|t|\leq X}
\Dist(f, n\mapsto \chi(n)n^{it};X)^2,$$
where the first infimum is taken over all Dirichlet characters $\chi$ of modulus $q\leq Q_{X,H'}$.

\begin{theorem}[{\cite[Theorem 1.7]{MRT15}}]\label{MRT_thm}
Let $X\geq H'\geq 10$ and $f$ be a $1$-bounded multiplicative function. Then we have
$$
    \sup_{\alpha\in[0,1]}\frac{1}{X}\int_0^X \Big|\Esub{j\in[H']} f(x+j)e(j\alpha)\Big| \d x
    \ll \exp\bigg(-\frac{M(f;X,Q_{X,H'})}{20}\bigg) + \frac{\log\log{H'}}{\log{H'}} + \frac{1}{(\log{X})^\frac{1}{700}}.
$$
\end{theorem}

In view of \cref{Cesaro_to_log}, we may apply this result with $\frac{1}{X}\int_0^X$ replaced by a logarithmic average.

We are ready to give the proof of \cref{large_xi_thm2}.

\begin{proof}[Proof of \cref{large_xi_thm2}]
We start by expanding the square to get
\begin{multline}\label{expanding_square_eq1}
    \sum_{\substack{\xi\in I_N \\ |\xi|>B}}
    \logE_{n\in[N]}\Big|\logE_{p\in\Pc_N}\f(n+p)\Big|^2 
    = \sum_{\substack{\xi\in I_N \\ |\xi|>B}}
    \logE_{n\in[N]}\logE_{p\in\Pc_N}h_{\xi,N}(n)\f(n+p)  \\
    = \sum_{\substack{\xi\in I_N \\ |\xi|>B}}
    \logE_{n\in[N]}\logE_{p\in\Pc_N}\f(n)
    \overline{h_{\xi,N}(n-p)}
    + \Oh\bigg(\frac{|I_N|\log{H}}{\log{N}}\bigg), 
\end{multline}
where we set $h_{\xi,N}(n) = \logE_{q\in\Pc_N}\f(n+q)$.
For the main term, we first decompose the set $\Pc_N$ into pairwise disjoint dyadic intervals $\Pc_k=[\frac{H_k}{2},H_k)\cap\P$, where $H_k = 2^k H_0$, for $k=1,\ldots,K$, where $K=\left\lceil\frac{\log(H/H_0)}{\log{2}}\right\rceil$.
Therefore, we may write it as
\begin{equation}\label{eq_main_term1}
    \sum_{\substack{\xi\in I_N \\ |\xi|>B}}\E^\ast_{k\in[K]}\logE_{n\in[N]}\logE_{p\in\Pc_k}\f(n)\overline{h_{\xi,N}(n-p)},
\end{equation}
where $\E^\ast := \frac{1}{\mathscr{L}}\sum_{k=1}^K \mathscr{L}_k$ and $\mathscr{L}_k = \sum_{p\in\Pc_k}\frac{1}{p} \sim \log 2$.
We define for $\xi\in I_N$, $|\xi|>B$,
$$\epsilon_\xi =
\begin{cases}
    \exp\Big(-\frac{\pi^2\xi^2}{200A^2}\Big), & |\xi|\leq\sqrt{|I_N|}, \\
    \Big(\frac{\log\log N}{\log\log\log N}\Big)^\frac{1}{5}\exp\Big(-\frac{3\log\log N}{20\log\log\log N}\Big), & \text{otherwise}.
\end{cases}$$
For simplicity, we denote the second value of $\epsilon_\xi$ by $\epsilon_N$.
For each $\xi,k$, applying \cref{bound_dyadic_set} for $f=\f$, $g=h_{\xi,N}$, $\xi\in I_N$, $\Pc=\Pc_k$ and $\epsilon=\epsilon_\xi>0$, it follows that \eqref{eq_main_term1} is

\begin{multline*}
    \ll \E^\ast_{k\in[K]}\sum_{\substack{\xi\in I_N \\ |\xi|>B}}\Bigg(\epsilon_\xi^{-4}\bigg(\sup_{\alpha\in[0,1]} \logE_{n\in[N]}\Big|\E_{j\in[H_k]}\f(n+j)e(j\alpha)\Big|+\frac{\log{H_k}}{\log N}\bigg)+\epsilon_\xi\Bigg) \\
    \ll \sum_{\substack{\xi\in I_N \\ |\xi|>B}}\epsilon_\xi^{-4} \E^\ast_{k\in[K]}\sup_{\alpha\in[0,1]}\logE_{n\in[N]}\Big|\Esub{j\in[H_k]} \f(n+j)e(j\alpha)\Big| + \frac{\log H}{\log N}\sum_{\substack{\xi\in I_N \\ |\xi|>B}}\epsilon_\xi^{-4} + \sum_{\substack{\xi\in I_N \\ |\xi|>B}}\epsilon_\xi.
\end{multline*}
Now, by \cref{MRT_thm}, 
\begin{multline*}
    \E^\ast_{k\in[K]}\sup_{\alpha\in[0,1]}\logE_{n\in[N]}\Big|\Esub{j\in[H_k]} \f(n+j)e(j\alpha)\Big| \\
    \ll \sup_{k\in[K]} \exp\bigg(-\frac{M(\f;N,Q_{N,H_k})}{20}\bigg) + \frac{\log\log H}{\log H_0} + \frac{1}{(\log N)^\frac{1}{700}}.
\end{multline*}
Combining Lemmas \ref{dist} and \ref{dist with chi}, summing over $\xi$, and noting that $\frac{\log\log H}{\log H_0}$ dominates over $\frac{\log H}{\log N}$ and $(\log N)^{-\frac{1}{700}}$, we obtain that \eqref{eq_main_term1} is

\begin{align*}
    \ll \sum_{B<\xi\leq\sqrt{|I_N|}} \epsilon_\xi^{-4}\exp\bigg(-\frac{\pi^2\xi^2}{40A^2}\bigg) + (\log N)^{-\frac{\pi^2}{A}} \sum_{\sqrt{|I_N|}<\xi\leq\frac{|I_N|}{2}} \epsilon_\xi^{-4}
    + \frac{\log\log H}{\log H_0}\sum_{B<\xi\leq \frac{|I_N|}{2}}\epsilon_\xi^{-4} & \\
    + \sum_{B<\xi\leq\frac{|I_N|}{2}}\epsilon_\xi & \\
    \ll \sum_{\xi>B}\exp\bigg(-\frac{\pi^2\xi^2}{200A^2}\bigg) + (\log N)^{-\frac{\pi^2}{A}}|I_N|\epsilon_N^{-4} + \frac{\log\log H}{\log H_0} \sum_{\xi\leq\sqrt{|I_N|}}\exp\bigg(\frac{\pi^2\xi^2}{50A^2}\bigg) & \\
    + \frac{\log\log H}{\log H_0}|I_N|\epsilon_N^{-4} 
    + |I_N|\epsilon_N. &
\end{align*}
Since $\frac{\log\log H}{\log H_0}\gg(\log N)^{-\frac{\pi^2}{A}}$, the second term is negligible. The last two terms are balanced by the choice of $\epsilon_N$ and they much smaller than $(\log\log N)^{-\Oh(1)}$, so they are also negligible. The first term, as in the proof of \cref{large_xi_thm1}, is 
$$\exp\bigg(-\frac{\pi^2B^2}{200A^2}\bigg) = \frac{1}{\log\log N}.$$
Finally, by a similar argument, we have the asymptotic
$$\sum_{\xi\leq\sqrt{|I_N|}}\exp\bigg(\frac{\pi^2\xi^2}{50A^2}\bigg) \sim \frac{25A^2}{\pi^2\sqrt{|I_N|}}\exp\bigg(\frac{\pi^2|I_N|}{50A^2}\bigg) \ll \exp\bigg(\frac{|I_N|}{2A^2}\bigg) = \exp\bigg(\sqrt{\frac{\log\log N}{\log\log\log N}}\bigg),$$
hence the third term is 
$$\ll \frac{\log\log N}{\log\log\log N}\exp\bigg(-\frac{3\log\log N}{4\log\log\log N} + \sqrt{\frac{\log\log N}{\log\log\log N}}\bigg)
\ll (\log\log N)^{-\Oh(1)}.$$
Putting everything together and since the error term in \eqref{expanding_square_eq1} is dominated by the main term, we conclude the proof.
\end{proof}

%SECTION
\section{Equivalence of Theorems~\ref{thm_A} and~\ref{thm_C}}
\label{Proof_Thm_C}

In this section we prove Theorem~\ref{thm_C}. In fact, we show that Theorems~\ref{thm_A} and~\ref{thm_C} are equivalent.

For our purposes in this section, we reformulate \eqref{eqn1_thm_C} in Theorem~\ref{thm_C} as follows: Given a bounded function $a\colon\N\to\C$, in light of \eqref{pk_and_logpk}, the left-hand side of \eqref{eqn1_thm_C} is 
\begin{align*}
    & \sum_{\ell=1}^\infty\overline{\pi}_\ell(N)\Big|\logE_{n\in[N]\cap\P_\ell}a(\Omega(n+1)) - \Esub{n\in[N]}a(\Omega(n))\Big| \\
    & \hspace*{0.2cm} = 
    \sum_{\ell=1}^\infty\Big|\Big(\logE_{n\in[N]\cap\P_\ell}a(\Omega(n+1))\Big)\overline{\pi}^{\log}_\ell(N) - \Big(\Esub{n\in[N]}a(\Omega(n))\Big)\overline{\pi}_\ell(N)\Big| + \Oh\Bigg(\frac{(\log\log\log{N})^2}{\sqrt{\log\log{N}}}\Bigg).
\end{align*}
Now observe that 
$$\Big(\logE_{n\in[N]\cap\P_\ell}a(\Omega(n+1))\Big)\overline{\pi}^{\log}_\ell(N)
= \logE_{n\in[N]} a(\Omega(n+1))\1_{\P_\ell}(n),$$
hence \eqref{eqn1_thm_C} is equivalent to
\begin{equation}\label{eqn_thm_C_reformulation}
    \sum_{\ell=1}^\infty \Big|\logE_{n\in[N]}a(\Omega(n+1))\1_{\P_\ell}(n) - \Big(\Esub{n\in[N]}a(\Omega(n))\Big)\overline{\pi}_\ell(N)\Big|
    \ll \frac{(\log\log\log{N})^2}{\sqrt{\log\log{N}}}.
\end{equation}

We first show the easy direction that Theorem~\ref{thm_C} implies Theorem~\ref{thm_A}.

\begin{proof}[Proof of Theorem~\ref{thm_C}$\implies$Theorem~\ref{thm_A}]
    Let $a,b\colon\N\to\C$ be bounded and suppose that Theorem~\ref{thm_C} is true. Then, for $N$ sufficiently large, applying \eqref{eqn_thm_C_reformulation} to $b$, we have
    \begin{align*}
        \logE_{n\in[N]} a(\Omega(n))b(\Omega(n+1))
        & = \sum_{\ell=1}^\infty a(\ell)\Big(\logE_{n\in[N]}\1_{\P_\ell}(n)b(\Omega(n+1))\Big) \\
        & = \sum_{\ell=1}^\infty a(\ell)\Big(\logE_{n\in[N]}\1_{\P_\ell}(n)b(\Omega(n+1))\Big) + \Oh\bigg(\frac{(\log\log\log{N})^2}{\sqrt{\log\log{N}}}\bigg) \\
        & = \Bigg(\sum_{\ell=1}^\infty a(\ell)\overline{\pi}^{\log}_{\ell}(N)\Bigg) \Big(\logE_{n\in[N]}b(\Omega(n))\Big) + \Oh\bigg(\frac{(\log\log\log{N})^2}{\sqrt{\log\log{N}}}\bigg) \\
        & = \Big(\logE_{n\in[N]}a(\Omega(n))\Big) \Big(\logE_{n\in[N]}b(\Omega(n))\Big) + \Oh\bigg(\frac{(\log\log\log{N})^2}{\sqrt{\log\log{N}}}\bigg), 
    \end{align*}
    where the last equality follows from \eqref{eqn_expanding_a(Omega(n))_2}.
    Finally, in light of \cref{aOmega_and_logaOmega}, we can substitute both logarithmic averages in the last expression above by the corresponding \Cesaro{} averages. This concludes the proof.
\end{proof}

Finally, we show the inverse implication, thus, establishing Theorem~\ref{thm_C}.

\begin{proof}[Proof of Theorem~\ref{thm_A}$\implies$Theorem~\ref{thm_C}]
Let $N\in\N$ be sufficiently large and $a\colon\N\to\C$. We can assume that $a$ is $1$-bounded.
We prove \eqref{eqn1_thm_C} and then \eqref{eqn2_thm_C} follows by Markov's inequality combined with the quantitative \Erdos-Kac theorem \cite[Theorem 2]{renyi-turan}. As we saw, it suffices to prove \eqref{eqn_thm_C_reformulation}.

By Theorem~\ref{thm_A}, we have that for any $1$-bounded function $b\colon\N\to\C$,
$$\Big|\logE_{n\in[N]}a(\Omega(n+1))b(\Omega(n)) - \Big(\Esub{n\in[N]}a(\Omega(n))\Big)\Big(\Esub{n\in[N]}b(\Omega(n)\Big)\Big| 
\ll \frac{(\log\log\log{N})^2}{\sqrt{\log\log{N}}}.$$

Therefore it is enough to prove the following:
If there is $\epsilon>0$ and $N_1<N_2<\ldots\in\N$, such that for all $i\in\N$,
\begin{equation}\label{assumption_in_proof_of_thmC}
    \sum_{\ell=1}^\infty \Big|\logE_{n\in[N_i]} a(\Omega(n+1))\1_{\P_\ell}(n) - \Big(\Esub{n\in[N_i]}a(\Omega(n))\Big)\overline{\pi}_\ell(N_i)\Big|
    \geq \epsilon,
\end{equation}
then there exists some $1$-bounded function $b\colon\N\to\C$, such that for all $i\in\N$,
\begin{equation}\label{conclusion_in_proof_of_thmC}
    \Big|\logE_{n\in[N_i]}a(\Omega(n+1))b(\Omega(n)) - \Big(\Esub{n\in[N_i]}a(\Omega(n))\Big)\Big(\Esub{n\in[N_i]}b(\Omega(n)\Big)\Big| \geq \epsilon.
\end{equation}
Observe that for $\ell>N_i$, $\overline{\pi}_\ell(N_i) = 0$ and $\1_{\P_\ell}(n) = 0$ for all $n$, hence \eqref{assumption_in_proof_of_thmC} becomes
\begin{equation}
\label{assumption2_in_proof_of_thmC}
    \sum_{\ell=1}^{N_i} \Big|\logE_{n\in[N_i]} a(\Omega(n+1))\1_{\P_\ell}(n) - \Big(\Esub{n\in[N_i]}a(\Omega(n))\Big)\overline{\pi}_\ell(N_i)\Big|
    \geq \epsilon.
\end{equation}
By passing to a subsequence of $(N_i)_{i\in\N}$ if necessary, we can assume without loss of generality that
\begin{align}
\label{eqn_A_implies_first_half_of_Aprimeprimeprime_2}
\sum_{\ell=1}^{N_{i-1}}\1_{\P_\ell}(n) \leq \frac{\epsilon}{4}
\quad\text{and}\quad
\sum_{\ell=1}^{N_{i-1}} \overline{\pi}_\ell(N_i)\leq \frac{\epsilon}{4}.
\end{align}
Let $\operatorname{sgn}\colon \R\to \{-1,0,1\}$ denote the function
\[
\operatorname{sgn}(x)=
\begin{cases}
-1,&\text{if}~x<0,\\
0,&\text{if}~x=0,\\
1,&\text{if}~x>0,\\
\end{cases}
\]
and we define $b\colon\N \to \{-1,0,1\}$ as
\[
b(\ell)=\operatorname{sgn}\Big(\logE_{n\in[N_i]} a(\Omega(n+1))\1_{\P_\ell}(n) - \Big(\Esub{n\in[N_i]}a(\Omega(n))\Big)\overline{\pi}_\ell(N_i)\Big)
\]
if $\ell\in (N_{i-1},N_i]$ for some $i\geq 2$ and $b(\ell)=0$ otherwise.

First, we invoke \eqref{eqn_A_implies_first_half_of_Aprimeprimeprime_2} and the triangle inequality to deduce that
\begin{multline*}
\sum_{\ell=1}^{N_{i-1}} \Big|\logE_{n\in[N_i]}  a(\Omega(n+1))\1_{\P_\ell}(n) -\Big(\Esub{n\in[N_i]}a(\Omega(n))\Big)\overline{\pi}_\ell(N_i)\Big| \\
\leq 
\Big|\logE_{n\in[N_i]}a(\Omega(n+1))\Big|\Bigg(\sum_{\ell=1}^{N_{i-1}} \1_{\P_\ell}(n)\Bigg) + \Big|\Esub{n\in[N_i]}a(\Omega(n))\Big|\Bigg(\sum_{\ell=1}^{N_{i-1}} \overline{\pi}_\ell(N_i)\Bigg) 
\leq
\frac{\epsilon}{2}.
\end{multline*}
Combined with \eqref{assumption2_in_proof_of_thmC}, this gives
\[
\sum_{\ell=N_{i-1}+1}^{N_i} \Big|\logE_{n\in[N_i]} a(\Omega(n+1))\1_{\P_\ell}(n)
- \Big(\Esub{n\in[N_i]} a(\Omega(n))\Big)\overline{\pi}_\ell(N_i)\Big|
\geq \frac{\epsilon}{2}.
\]
Using the definition of $b(\ell)$, it now follows that 
\[
\sum_{\ell=N_{i-1}+1}^{N_i} b(\ell)\bigg(\logE_{n\in[N_i]} a(\Omega(n+1))\1_{\P_\ell}(n) 
- \Big(\Esub{n\in[N_i]} a(\Omega(n))\Big)\overline{\pi}_\ell(N_i)\bigg)
\geq \frac{\epsilon}{2}.
\]
Using \eqref{eqn_A_implies_first_half_of_Aprimeprimeprime_2} 
and the triangle inequality again we get
\begin{multline}\label{last_step_pf_thm_C}
    \Bigg|\sum_{\ell=1}^\infty b(\ell) \bigg(\logE_{n\in[N_i]} a(\Omega(n+1))\1_{\P_\ell}(n) - \Big(\Esub{n\in[N_i]} a(\Omega(n))\Big)\overline{\pi}_\ell(N_i)\bigg)\Bigg| \\
    = \Bigg|\sum_{\ell=1}^{N_i} b(\ell) \bigg(\logE_{n\in[N_i]} a(\Omega(n+1))\1_{\P_\ell}(n) - \Big(\Esub{n\in[N_i]} a(\Omega(n))\Big)\overline{\pi}_\ell(N_i)\bigg)\Bigg|
    \geq \epsilon.
\end{multline}
By rearranging the first expression above and then invoking \eqref{eqn_expanding_a(Omega(n))_1} and \eqref{eqn_expanding_a(Omega(n))_2}, 
we have 
\begin{align*}
    \sum_{\ell=1}^\infty b(\ell) & \bigg(\logE_{n\in[N_i]} a(\Omega(n+1))\1_{\P_\ell}(n) - \Big(\Esub{n\in[N_i]} a(\Omega(n))\Big)\overline{\pi}_\ell(N_i)\bigg) \\
    & = \logE_{n\in[N_i]}a(\Omega(n+1))\Bigg(\sum_{\ell=1}^\infty b(\ell)\1_{\P_\ell}(n)\Bigg) - \Big(\Esub{n\in[N_i]}a(\Omega(n))\Big)\Bigg(\sum_{\ell=1}^\infty b(\ell)\overline{\pi}_\ell(N_i)\Bigg) \\
    & = \logE_{n\in[N_i]}a(\Omega(n+1))b(\Omega(n)) - \Big(\Esub{n\in[N_i]}a(\Omega(n))\Big)\Big(\Esub{n\in[N_i]}b(\Omega(n))\Big).
\end{align*}
Combining this with \eqref{last_step_pf_thm_C}, \eqref{conclusion_in_proof_of_thmC} follows, thus, concluding the proof.
\end{proof}

\appendix
\section{Another proof of Theorem~\ref{thm_A} with optimal error term}
\label{appendix}

We sketch an alternative proof of Theorem~\ref{thm_A} with optimal error term that was communicated to us by Harald A.~Helfgott after the first version of the article appeared. The three main ingredients of this proof is the main theorem of \cite{helf_radz} (in particular, \cite[Corollary 1.4]{helf_radz}), the quantitative \Erdos-Kac theorem \cite[Theorem 2]{renyi-turan}, and \cite[Proposition 8.12]{helf_radz}. The use of the latter as a black-box is crucial to obtain the optimal error term. 

Let $a,b\colon\N\to\C$ be bounded functions. We want to show 
\begin{equation*}
    \Big|\logE_{n\in[N]}a(\Omega(n))b(\Omega(n+1)) -\Big(\E_{n\in[N]}a(\Omega(n))\Big)\Big(\E_{n\in[N]}b(\Omega(n))\Big)\Big|
    \ll \frac{1}{\sqrt{\log\log N}}.
\end{equation*}
In view of \cref{aOmega_and_logaOmega}, it is enough to show
\begin{equation*}
    \Big|\logE_{n\in[N]}a(\Omega(n))b(\Omega(n+1)) -\Big(\logE_{n\in[N]}a(\Omega(n))\Big)\Big(\logE_{n\in[N]}b(\Omega(n))\Big)\Big|
    \ll \frac{1}{\sqrt{\log\log N}}.
\end{equation*}

We let $H = \exp\big(\sqrt{\log N}\big)$, $H_0 = \exp\big((\log N)^\frac{4}{9}\big)$, $\Pc = [H_0,H]\cap\P$, $\mathscr{L} = \sum_{p\in\Pc}\frac{1}{p} \asymp \log\log N$, and $T_N = (\log\log{N} - A\sqrt{\log\log{N}}, \log\log{N} + A\sqrt{\log\log{N}}]\cap\Z$, with $A=\sqrt{\log\log\log N}$ as in \eqref{I_N_def}.

We start by noting that 
\begin{equation*}
    \logE_{n\in[N]}a(\Omega(n))b(\Omega(n+1)) = \logE_{p\in\Pc}\logE_{n\in[N]}\widetilde{a}(\Omega(n))\widetilde{b}(\Omega(n+p)) + \Oh\bigg(\frac{1}{\sqrt{\log\log N}}\bigg),
\end{equation*}
where $\widetilde{a},\widetilde{b}$ are defined as in the proof of \cref{reducing_thm_A}.
The proof of this is identical to that of \cref{Helf-Radz_cor}, and is a consequence of \cite[Corollary 1.4]{helf_radz}.
Then it suffices to show 
\begin{equation*}
    \Big|\logE_{p\in\Pc}\logE_{n\in[N]}\widetilde{a}(\Omega(n))\widetilde{b}(\Omega(n+p)) - \Big(\logE_{n\in[N]}\widetilde{a}(\Omega(n))\Big)\Big(\logE_{n\in[N]}\widetilde{b}(\Omega(n))\Big)\Big| \ll \frac{1}{\sqrt{\log\log N}}.
\end{equation*}
We can write the left-hand side of the last equation as
\begin{multline*}
    \Bigg|\sum_{\ell,k=1}^\infty \widetilde{a}(\ell)\widetilde{b}(k)\logE_{p\in\Pc}\logE_{n\in[N]}\1_{\P_\ell}(n)\1_{\P_k}(n+p) \\
    - \Bigg(\sum_{\ell=1}^\infty \widetilde{a}(\ell)\logE_{n\in[N]}\1_{\P_\ell}(n)\Bigg)\Bigg(\sum_{k=1}^\infty \widetilde{b}(k)\logE_{n\in[N]}\1_{\P_k}(n)\Bigg)\Bigg|,
\end{multline*}
hence it suffices to show 
\begin{equation*}
    \sum_{\ell,k=1}^\infty \Big|\logE_{p\in\Pc}\logE_{n\in[N]}\1_{\P_\ell}(n)\1_{\P_k}(n+p) - \Big(\logE_{n\in[N]}\1_{\P_\ell}(n)\Big)\Big(\logE_{n\in[N]}\1_{\P_k}(n)\Big)\Big| \ll \frac{1}{\sqrt{\log\log N}}.
\end{equation*}
Now we restrict both summations to $\ell,k\in T_N$ and the contribution of the complementary range is at most
$$2\logE_{n\in[N]} \1_{\Omega(n)\not\in T_N} + \logE_{p\in\Pc}\logE_{n\in[N]}\1_{\Omega(n+p)\not\in T_N}
\ll \E_{n\in[N]}\1_{\Omega(n)\not\in T_N} + \sup_{p\in\Pc}\E_{n\in[N]}\1_{\Omega(n+p)\not\in T_N}.$$
By \cref{renyi-turan_cor}, which is a corollary of the quantitative \Erdos-Kac theorem \cite[Theorem 2]{renyi-turan}, the error is
$$\ll \frac{1}{\sqrt{\log\log N}},$$
thus, it is enough to show 
\begin{equation}\label{app_eq}
    \sum_{\ell,k\in T_N}\Big|\logE_{p\in\Pc}\logE_{n\in[N]}\1_{\P_\ell}(n)\1_{\P_k}(n+p) - \Big(\logE_{n\in[N]}\1_{\P_\ell}(n)\Big)\Big(\logE_{n\in[N]}\1_{\P_k}(n)\Big)\Big| \ll \frac{1}{\sqrt{\log\log N}}.
\end{equation}
By \cite[Proposition 8.12]{helf_radz}, we have for any $M\in[N^\frac{1}{\log\log N},N^{2\log\log N}]$ and uniformly for all $k,\ell\geq1$,
\begin{multline*}
    \logE_{n\in\Pc}\Esub{n\in[M]}\1_{\P_\ell}(n)\1_{\P_k}(n+p) \\
    = \Big(\Esub{n\in[M]}\1_{\P_\ell}(n)\Big)\Big(\Esub{n\in[M]}\1_{\P_k}(n)\Big)
    + \Oh\Bigg(\frac{|\ell-\log\log N|^2+|k-\log\log N|^2+1}{(\log\log N)^2}\Bigg).
\end{multline*}
Using \cref{Cesaro_to_log} as we did in \cref{small_xi_section} to derive \cref{small_xi_thm} from \cref{small_xi_Cesaro_thm}, it follows that 
\begin{multline*}
    \logE_{n\in\Pc}\logE_{n\in[N]}\1_{\P_\ell}(n)\1_{\P_k}(n+p) \\
    = \Big(\logE_{n\in[N]}\1_{\P_\ell}(n)\Big)\Big(\logE_{n\in[N]}\1_{\P_k}(n)\Big)
    + \Oh\Bigg(\frac{|\ell-\log\log N|^2+|k-\log\log N|^2+1}{(\log\log N)^2}\Bigg).
\end{multline*}
Summing over all $\ell,k\in T_N$, we obtain that the left-hand side of \eqref{app_eq} is 
$$\ll \frac{\log\log\log N}{\log\log N},$$
so the result follows.

While an optimal error term for the \Erdos{}-Kac theorem was obtained by R\'enyi and \Turan{} in \cite{renyi-turan}, it has remained an open problem to do the same for the bivariate \Erdos{}-Kac theorem.
Currently, the best known result, established in \cite{Mangerel20}, asserts that
\begin{multline*}
        \E_{n\in[N]}\1_{(-\infty,y)}\bigg(\frac{\Omega(n)-\log\log N}{\sqrt{\log\log N}}\bigg)\1_{(-\infty,z)}\bigg(\frac{\Omega(n+1)-\log\log N}{\sqrt{\log\log N}}\bigg) \\
        = \Bigg(\frac{1}{\sqrt{2\pi}}\int_{-\infty}^y e^{-\frac{t^2}{2}}\d t\Bigg)\Bigg(\frac{1}{\sqrt{2\pi}}\int_{-\infty}^z e^{-\frac{t^2}{2}}\d t\Bigg) + \Oh\Bigg(\frac{\log\log\log N}{(\log\log N)^\frac{1}{4}}\Bigg).
\end{multline*}
As a consequence of the above argument, we can now establish an optimal error in the bivariate setting, but with the caveat that \Cesaro{} averages are replaced by logarithmic averages.

\begin{corollary}
    Uniformly for all $y,z\in\R$ and for $N$ sufficiently large, we have
    \begin{multline*}
        \logE_{n\in[N]}\1_{(-\infty,y)}\bigg(\frac{\Omega(n)-\log\log N}{\sqrt{\log\log N}}\bigg)\1_{(-\infty,z)}\bigg(\frac{\Omega(n+1)-\log\log N}{\sqrt{\log\log N}}\bigg) \\
        = \Bigg(\frac{1}{\sqrt{2\pi}}\int_{-\infty}^y e^{-\frac{t^2}{2}}\d t\Bigg)\Bigg(\frac{1}{\sqrt{2\pi}}\int_{-\infty}^z e^{-\frac{t^2}{2}}\d t\Bigg) + \Oh\bigg(\frac{1}{\sqrt{\log\log N}}\bigg).
    \end{multline*}
\end{corollary}

\bibliographystyle{amsalpha}
\bibliography{Refs_Omega.bib}

@Article{Loyd23,
  author   = {Loyd, K.},
  title    = {{A} dynamical approach to the asymptotic behavior of the sequence {$\Omega(n)$}},
  journal  = {Ergodic Theory Dynam. Systems},
  year     = {2023},
  volume   = {43},
  number   = {11},
  pages    = {3685--3706},
  fjournal = {Ergodic Theory and Dynamical Systems}
}

@article{Hal,
  title={{\"U}ber die {M}ittelwerte multiplikativer zahlentheoretischer {F}unktionen},
  author={Hal{\'a}sz, G.},
  journal={Acta Mathematica Hungarica},
  volume={19},
  number={3-4},
  pages={365--403},
  year={1968},
  publisher={Akad{\'e}miai Kiad{\'o}, co-published with Springer Science+ Business Media BV~…}
}

@article{richert,
 AUTHOR = {Richert, H.-E.},
     TITLE = {Zur {A}bsch\"atzung der {R}iemannschen {Z}etafunktion in der
              {N}\"ahe der {V}ertikalen {$\sigma =1$}},
   JOURNAL = {Math. Ann.},
  FJOURNAL = {Mathematische Annalen},
    VOLUME = {169},
      YEAR = {1967},
     PAGES = {97--101}
}

@Article{MRT15,
  author   = {Matom\"{a}ki, K. and Radziwi{\l}{\l}, M. and Tao, T.},
  title    = {An averaged form of {C}howla's conjecture},
  journal  = {Algebra Number Theory},
  year     = {2015},
  volume   = {9},
  number   = {9},
  pages    = {2167--2196},
  fjournal = {Algebra \& Number Theory}
}

@Article{Tao16,
  author   = {Tao, T.},
  title    = {The logarithmically averaged {C}howla and {E}lliott conjectures for two-point correlations},
  journal  = {Forum Math. Pi},
  year     = {2016},
  volume   = {4},
  pages    = {e8, 36},
  fjournal = {Forum of Mathematics. Pi}
}

@book {Montogomery_ten,
    AUTHOR = {Montgomery, H. L.},
     TITLE = {Ten lectures on the interface between analytic number theory
              and harmonic analysis},
    SERIES = {CBMS Regional Conference Series in Mathematics},
    VOLUME = {84},
 PUBLISHER = {Conference Board of the Mathematical Sciences, Washington, DC;
              by the American Mathematical Society, Providence, RI},
      YEAR = {1994},
}

@book {Tenenbaum_book,
    AUTHOR = {Tenenbaum, G.},
     TITLE = {Introduction to analytic and probabilistic number theory},
    SERIES = {Graduate Studies in Mathematics},
    VOLUME = {163},
   EDITION = {Third},
    PUBLISHER = {American Mathematical Society, Providence, RI},
      YEAR = {2015}
}

@Article{TT18,
  author   = {Tao, T. and Ter\"{a}v\"{a}inen, J.},
  title    = {Odd order cases of the logarithmically averaged {C}howla conjecture},
  journal  = {J. Th\'{e}or. Nombres Bordeaux},
  year     = {2018},
  volume   = {30},
  number   = {3},
  pages    = {997--1015}
}

@article{Helf_Ubis,
    AUTHOR = {Helfgott, H. A. and Ubis, A.},
     TITLE = {Primes, parity and analysis},
   JOURNAL = {Publ. Mat. Urug.},
  FJOURNAL = {Publicaciones Matem\'aticas del Uruguay},
    VOLUME = {18},
      YEAR = {2023},
     PAGES = {205--283}
}

@article{helf_radz,
  title={Expansion, divisibility and parity},
  author={Helfgott, H. A. and Radziwi{\l}{\l}, M.},
  journal={arXiv preprint arXiv:2103.06853},
  year = {2021}
}

@article{renyi-turan,
    author     = {R\'{e}nyi, A. and Tur\'{a}n, P.},
  title      = {On a theorem of {E}rd{\H{o}}s-{K}ac},
  journal    = {Acta Arith.},
  year       = {1958},
  volume     = {4},
  pages      = {71--84}
}

@book{MV_book,
    AUTHOR = {Montgomery, H. L. and Vaughan, R. C.},
     TITLE = {Multiplicative number theory. {I}. {C}lassical theory},
    SERIES = {Cambridge Studies in Advanced Mathematics},
    VOLUME = {97},
 PUBLISHER = {Cambridge University Press, Cambridge},
      YEAR = {2007}
}

@article{elliott_TK,
 AUTHOR = {Elliott, P. D. T. A.},
     TITLE = {The {T}ur\'an-{K}ubilius inequality},
   JOURNAL = {Proc. Amer. Math. Soc.},
  FJOURNAL = {Proceedings of the American Mathematical Society},
    VOLUME = {65},
      YEAR = {1977},
    NUMBER = {1},
     PAGES = {8--10}
}

@article{pilatte2023improved,
  title={Improved bounds for the two-point logarithmic {C}howla conjecture},
  author={Pilatte, C.},
  journal={arXiv preprint arXiv:2310.19357},
  year={2023}
}

@incollection{tao2017equivalence,
      AUTHOR = {Tao, T.},
     TITLE = {Equivalence of the logarithmically averaged {C}howla and
              {S}arnak conjectures},
 BOOKTITLE = {Number theory--{D}iophantine problems, uniform distribution
              and applications},
     PAGES = {391--421},
 PUBLISHER = {Springer, Cham},
      YEAR = {2017}
}

@article{MSTT2023,
     AUTHOR = {Matom\"aki, K. and Shao, X. and Tao, T. and
              Ter\"av\"ainen, J.},
     TITLE = {Higher uniformity of arithmetic functions in short intervals
              {I}. {A}ll intervals},
   JOURNAL = {Forum Math. Pi},
  FJOURNAL = {Forum of Mathematics. Pi},
    VOLUME = {11},
      YEAR = {2023},
     PAGES = {Paper No. e29, 97}
}

@article{Matomaki-Radziwill-Shao-Tao-Teravainen2026,
    author  = {K. Matom{\"a}ki and M. Radziwi{\l}{\l} and X. Shao and T. Tao and J. Ter{\"a}v{\"a}inen},
    title   = {Higher uniformity of arithmetic functions in short intervals {II}. Almost all intervals},
    journal = {Invent. Math.},
    volume  = {244},
    number  = {3},
    pages   = {967--1091},
    year    = {2026}
}

@article{KLRT2024,
    AUTHOR = {Kanigowski, A. and Lema\'nczyk, M. and Richter, F.
              K. and Ter\"av\"ainen, J.},
     TITLE = {On the local {F}ourier uniformity problem for small sets},
   JOURNAL = {Int. Math. Res. Not. IMRN},
  FJOURNAL = {International Mathematics Research Notices. IMRN},
      YEAR = {2024},
    NUMBER = {15},
     PAGES = {11488--11512}
}

@article{MRTT2023_higher_uniformity,
  AUTHOR = {Matom\"aki, K. and Radziwi\l\l, M. and Tao, T. and
              Ter\"av\"ainen, J. and Ziegler, T.},
     TITLE = {Higher uniformity of bounded multiplicative functions in short
              intervals on average},
   JOURNAL = {Ann. of Math. (2)},
  FJOURNAL = {Annals of Mathematics. Second Series},
    VOLUME = {197},
      YEAR = {2023},
    NUMBER = {2},
     PAGES = {739--857}
}

@article {Selberge_Sathe,
    AUTHOR = {Selberg, A.},
     TITLE = {Note on a paper by {L}. {G}. {S}athe},
   JOURNAL = {J. Indian Math. Soc. (N.S.)},
  FJOURNAL = {The Journal of the Indian Mathematical Society. New Series},
    VOLUME = {18},
      YEAR = {1954},
     PAGES = {83--87}
}

@article {Korobov_zero-free,
    AUTHOR = {Korobov, N. M.},
     TITLE = {Estimates of trigonometric sums and their applications},
   JOURNAL = {Uspehi Mat. Nauk},
  FJOURNAL = {Akademija Nauk SSSR i Moskovskoe Matemati\v ceskoe Ob\v s\v
              cestvo. Uspehi Matemati\v ceskih Nauk},
    VOLUME = {13},
      YEAR = {1958},
    NUMBER = {4(82)},
     PAGES = {185--192}
}

@article{Vinogradov_zero-free,
    AUTHOR = {Vinogradov, I. M.},
     TITLE = {A new estimate of the function {$\zeta (1+it)$}},
   JOURNAL = {Izv. Akad. Nauk SSSR Ser. Mat.},
  FJOURNAL = {Izvestiya Akademii Nauk SSSR. Seriya Matematicheskaya},
    VOLUME = {22},
      YEAR = {1958},
     PAGES = {161--164}
}

@article{Charamaras_MVT_24,
  author        = {Charamaras, D.},
  title         = {Mean value theorems in multiplicative systems and joint ergodicity of additive and multiplicative actions},
  JOURNAL = {Trans. Amer. Math. Soc.},
  FJOURNAL = {Transactions of the American Mathematical Society},
    VOLUME = {378},
      YEAR = {2025},
    NUMBER = {3},
     PAGES = {1883–1937}
}

@article{Richter_PNT_21,
    AUTHOR = {Richter, F. K.},
    TITLE = {A new elementary proof of the {P}rime {N}umber {T}heorem},
    JOURNAL = {Bull. Lond. Math. Soc.},
    FJOURNAL = {Bulletin of the London Mathematical Society},
    VOLUME = {53},
    YEAR = {2021},
    NUMBER = {5},
    PAGES = {1365--1375}
}

@Article{BR22,
  author   = {Bergelson, V. and Richter, F. K.},
  title    = {{D}ynamical generalizations of the prime number theorem and disjointness of additive and multiplicative semigroup actions},
  journal  = {Duke Math. J.},
  year     = {2022},
  volume   = {171},
  number   = {15},
  pages    = {3133--3200},
  fjournal = {Duke Mathematical Journal}
}

@Article{EAKLdlR17,
  author   = {El Abdalaoui, E. H. and Ku{\l}aga-Przymus, J. and Lema{\'n}czyk, M. and de la Rue, T.},
  title    = {The {C}howla and the {S}arnak conjectures from ergodic theory point of view},
  journal  = {Discrete Contin. Dyn. Syst.},
  year     = {2017},
  volume   = {37},
  number   = {6},
  pages    = {2899--2944},
  fjournal = {Discrete and Continuous Dynamical Systems. Series A}
}

@Article{Erdos48a,
  author  = {Erd\H{o}s, P.},
  title   = {On the integers having exactly {$K$} prime factors},
  journal = {Ann. of Math. (2)},
  year    = {1948},
  volume  = {49},
  pages   = {53--66}
}

@InCollection{FKL18,
  author    = {Ferenczi, S. and Ku{\l}aga-Przymus, J. and Lema{\'n}czyk, M.},
  title     = {Sarnak's conjecture: what's new},
  booktitle = {Ergodic theory and dynamical systems in their interactions with arithmetics and combinatorics},
  publisher = {Springer, Cham},
  year      = {2018},
  volume    = {2213},
  series    = {Lecture Notes in Math.},
  pages     = {163--235}
}

@Article{Frantzikinakis17,
  author   = {Frantzikinakis, N.},
  title    = {Ergodicity of the {L}iouville system implies the {C}howla conjecture},
  journal  = {Discrete Anal.},
  year     = {2017},
  pages    = {Paper No. 19, 41},
  fjournal = {Discrete Analysis}
}

@Article{FH18a,
  author   = {Frantzikinakis, N. and Host, B.},
  title    = {The logarithmic {S}arnak conjecture for ergodic weights},
  journal  = {Ann. of Math. (2)},
  year     = {2018},
  volume   = {187},
  number   = {3},
  pages    = {869--931}
}

@Article{GLdlR21,
  author        = {Gomilko, A. and Lema{\'n}czyk, M. and de la Rue, T.},
  title         = {On {F}urstenberg systems of aperiodic multiplicative functions of {M}atom\"{a}ki, {R}adziwi{\l}{\l} and {T}ao},
  journal       = {J. Mod. Dyn.},
  fjournal      = {Journal of Modern Dynamics},
  year          = {2021},
  volume        = {17},
  pages         = {529--555}
}

@InCollection{Ramare18,
  author    = {Ramar\'{e}, O.},
  title     = {Chowla's conjecture: from the {L}iouville function to the {M}oebius function},
  booktitle = {Ergodic theory and dynamical systems in their interactions with arithmetics and combinatorics},
  publisher = {Springer, Cham},
  year      = {2018},
  volume    = {2213},
  series    = {Lecture Notes in Math.},
  pages     = {317--323}
}

@Book{Chowla87,
  title     = {The {R}iemann {H}ypothesis and {H}ilbert's {T}enth {P}roblem},
  publisher = {CRC Press},
  year      = {1987},
  author    = {Chowla, S.},
  editor    = {Merrill's International Series in Engineering Technology},
}

@misc{Sarnak11,
    author       = {P. Sarnak},
    title        = {Three Lectures on the M{\"o}bius Function, Randomness and Dynamics},
    year         = {2009},
    note         = {IAS Lecture Notes},
    howpublished = {\url{https://publications.ias.edu/sites/default/files/MobiusFunctionsLectures%282%29_0.pdf}}
}

@Article{Elliott92,
  author   = {Elliott, P. D. T. A.},
  title    = {On the correlation of multiplicative functions},
  journal  = {Notas Soc. Mat. Chile},
  year     = {1992},
  volume   = {11},
  number   = {1},
  pages    = {1--11},
  issn     = {0716-1298},
  fjournal = {Notas de la Sociedad de Matem\'{a}tica de Chile},
}

@article{KMT23_almost_all_scales,
  title={On {E}lliott's conjecture and applications},
  author={Klurman, O. and Mangerel, A. P. and Ter{\"a}v{\"a}inen, J.},
  journal={arXiv preprint arXiv:2304.05344},
  year={2023}
}

@article{GS_large_char_sums,
    AUTHOR = {Granville, A. and Soundararajan, K.},
     TITLE = {Large character sums: {B}urgess's theorem and zeros of
              {$L$}-functions},
   JOURNAL = {J. Eur. Math. Soc. (JEMS)},
  FJOURNAL = {Journal of the European Mathematical Society (JEMS)},
    VOLUME = {20},
      YEAR = {2018},
    NUMBER = {1},
     PAGES = {1--14}}

@book{De_Bruijn_81,
  title={Asymptotic methods in analysis},
  author={De Bruijn, N. G.},
  volume={4},
  year={1981},
  publisher={Courier Corporation}
}

@article{Mangerel20,
  author  = {Mangerel, A. P.},
  title   = {On the bivariate Erd{\H{o}}s--Kac theorem and correlations of the M{\"o}bius function},
  journal = {Mathematical Proceedings of the Cambridge Philosophical Society},
  volume  = {169},
  number  = {3},
  year    = {2020},
  pages   = {547--605},
  doi     = {10.1017/S0305004119000288}
}

%AFFILIATION
\bigskip
\bigskip
\footnotesize
\noindent
Dimitrios Charamaras\\
\textsc{Centre de recherches math\'{e}matiques, Universit\'e de Montr\'eal}\\
\href{mailto:dimitrios.charamaras@epfl.ch}
{\texttt{dimitrios.charamaras@umontreal.ca}}

\bigskip
\noindent
Florian K. Richter\\
\textsc{{\'E}cole Polytechnique F{\'e}d{\'e}rale de Lausanne (EPFL)}\\
\href{mailto:f.richter@epfl.ch}
{\texttt{f.richter@epfl.ch}}

%END DOCUMENT
\end{document}